\documentclass{article}
\usepackage{graphicx}

\graphicspath{{./figures/}}
\usepackage{epstopdf}
\usepackage{amsmath,amssymb,amsfonts,subeqnarray}
\usepackage[usenames,dvipsnames]{color}
\usepackage{float}
\usepackage{subfig}
\usepackage{natbib}
\usepackage{todonotes}
\newcommand{\PIN}{\textnormal{{\tiny PIN}}}

\newcommand{\IAA}{\textnormal{{\tiny IAA}}}

\newcommand{\T}{\textnormal{{\tiny T}}}
\newcommand{\Tn}{\textnormal{T}}
\newcommand{\Dn}{\textnormal{D}}
\newcommand{\field}[1]{\mathbb{#1}} 
\newcommand{\R}{\field{R}} 
\providecommand{\change}{\textcolor{black}}
\begin{document}

%
%
%
%
\title{Numerical bifurcation analysis of the pattern formation in a cell based auxin transport model}


\author{Delphine Draelants \\
\small          Department of Mathemtics and Computer Science, Universiteit Antwerpen, \\
\small           Campus Middelheim, Building G, Middelheimlaan 1, B 2020, ANTWERPEN \\
\small           Tel.: +323-265-38-59 \\
\small           \tt{Delphine.Draelants@ua.ac.be}   
          \and
           Jan Broeckhove \\
\small          Department of Mathematics and Computer Science, Universiteit Antwerpen, \\
\small          Campus Middelheim, Building G, Middelheimlaan 1, B 2020, ANTWERPEN
          \and
           Gerrit T.S. Beemster \\
\small           Department of Biology, Universiteit Antwerpen, \\
\small           Campus Middelheim, Building U, Groenenborgerlaan 171, B 2020, ANTWERPEN
          \and
           Wim Vanroose \\
\small           Department of Mathematics and Computer Science, Universiteit Antwerpen, \\
\small           Campus Middelheim, Building G, Middelheimlaan 1, B 2020, ANTWERPEN
}

\date{Received: date / Accepted: date}

\maketitle

%
%
%
%
\begin{abstract}
  Transport models of growth hormones can be used to reproduce the
  hormone accumulations that occur in plant organs. Mostly, these
  accumulation patterns are calculated using time step methods, even
  though only the resulting steady state patterns of the model are of
  interest.  We examine the steady state solutions of the hormone
  transport model of \citet{Smith2006} for a one-dimensional row of
  plant cells. We search for the steady state solutions as a function
  of three of the model parameters by using numerical continuation
  methods and bifurcation analysis. These methods are more adequate
  for solving steady state problems than time step methods.  We
  discuss a trivial solution where the concentrations of hormones are
  equal in all cells and examine its stability region. We identify two
  generic bifurcation scenarios through which the trivial solution
  loses its stability.  The trivial solution becomes either a steady
  state pattern with regular spaced peaks or a pattern where the
  concentration is periodic in time.

  \textbf{keywords:}{Bifurcation analysis, pattern formation,
    parameter dependence, auxin transport model, stability, periodic
    solution pattern} 
\end{abstract}

%
%
%
%
\section{Introduction}
\label{sec:Introduction}

%
%
\subsection{Biological background}
\label{subsec:BiologicalBackground}
For centuries, the formation of well-defined patterns in plants, such
as the orientation and shape of leaves, their venation patterns, the
spatial distribution of hairs and stomata, the early embryonic
development patterns and the branching patterns in both root systems
and treetops, has intrigued many scientists. Experimental research has
identified a number of molecular components that play a major role in
several of these pattern formation processes. One of them is the plant
hormone auxin, and more specifically the auxin molecule Indol-3-Acetic
Acid (IAA). Experiments have shown that the active directional
transport, which leads to accumulation spots of the auxin hormone,
plays a central part in the pattern formation 
\citep{Scarpella2006,Benkova2003,Bilsborough2011}.  

Based on such experimental evidence,
Reinhardt~\textit{et~al.} developed a conceptual model that describes
the auxin transport through the cells \citep{Reinhardt2003}. Smith and
collaborators then constructed a computational simulation
model \citep{Smith2006} incorporating the experimental evidence that
the transport of the auxin molecule IAA is driven by a pumping
mechanism that is mediated by PIN1 proteins located at the cell
membrane in addition to diffusion \citep{Palme1999}.  Therefore, Smith
and collaborators modeled the transport of the IAA hormone through the
cells by describing the simultaneous evolution of the PIN1 protein and
the IAA hormone concentrations over time.  Also other computational
models were developed based on these molecular mechanisms identified
by Reinhardt~\textit{et~al.}.  For instance J\"{o}nsson \textit{et
  al.} proposed a phyllotaxis model based on the polarized auxin
transport \citep{Jonsson2006}.  They analyzed a simplified version of
their model that assumes an equal and constant PIN1 concentration 
in every cell and membrane.  In their simulations they used a linear 
row of uniform cells with periodic boundary conditions. The results 
show that the spacing and the number of peaks in this simplified model 
depends on the different parameters. J\"onsson and collaborators also 
performed a stability analysis and found an analytical expression for 
the eigenvalues, belonging to a solution pattern with equal auxin concentrations.
The eigenvalues are all real and a function of the model 
parameters of the model. They also identified the parameter threshold 
where the largest eigenvalue becomes unstable. Beyond this threshold,
all stable solutions will contain auxin peaks.

This paper expands the study of the steady states in the transport of
hormones. We limit ourselves to the study of the auxin distribution in
a linear row of uniform cells that represents, for example, a cross
section through a young leaf. 
We perform a thorough mathematical
exploration of the behavior of the models and how their equations are
solved starting from the basic coupled model of \citet{Smith2006}. In contrast to the analysis of
J\"onsson~\textit{et~al.} on a row of cells \change{with fixed concentration of PIN1}, 
we will use a coupled model where the PIN1 concentration is allowed
to change from cell to cell. The analysis gives new insights into the
spacing of auxin accumulations that form the basis of vascular
development \citep{Scarpella2006}.

\change{ The patterns that emerge in a dynamical system are often
  studied mathematically through a bifurcation analysis. It relates
  the stability of the patterns to the parameters that occur in the
  systems description.  The transitions where the patterns lose their stability are
  bifurcation points. For an overview  of the bifurcation analysis of
  patterns we refer to the book \citep{Hoyle}.}

The main contribution of the paper is a systematic numerical
bifurcation analysis for the coupled model of \citet{Smith2006} describing the transport of
\change{auxin}.  The analysis identifies two generic bifurcation
scenarios that reappear for various choices of the parameters of the
problem.  Through the bifurcation diagrams we identify the genesis of
the patterns that were observed by Smith and
collaborators. Furthermore, we have found a limited parameter range
that allows periodic solutions in the system. In these solutions, the
concentration of auxin of each cell varies periodically over time.  To
the knowledge of the authors these results have not appeared in the
literature.

We present our work as follows.  In section
\ref{subsec:DescriptionOfTheMathematicalModel} the basic cell polarization and auxin
transport model of \citet{Smith2006}, where PIN1 is allowed to
change, is reconstructed. \change{We also introduce a slightly generalized version
of the active transport equation of this model.} In the next section, a
specific model that will be used in the simulations is defined.  In
\citet{Smith2006} the domain roughly correlates to a ring of cells
around an axial plant organ. In subsequent situations we consider a
linear row of cells running from the margin to the midvein of the leaf
and we consider \change{zero} fluxes at the boundary of the leaf. We will
describe this by using \change{homogeneous} Neumann boundary conditions instead of periodic
boundary conditions which is explained in section
\ref{subsec:TheDomainItsBoundaryConditionsAndTheParameters}. Also in
this section we look at the different parameter values. Similar to
Smith \textit{et al.}, we use time integration to solve the coupled
equations of Smith \textit{et al.} in section
\ref{subsec:TimeIntegration}. Since we are only interested in the
steady state solutions, we define in section
\ref{subsec:SteadyStateProblem}, the corresponding steady state
\change{systems of the slightly different equations}.  For \change{these} steady state model\change{s}, we define in section 
\ref{sec:TheTrivialSolution} a trivial solution and its stability properties. 
The stability is dependent on the model parameters and we examine
for which parameter regions the trivial solution is stable. 
Section \ref{sec:Methods} contains the techniques that will
be used to solve the model\change{s}. In particular, we will discuss bifurcation
analysis (\ref{subsec:BifurcationAnalysis}) and continuation methods
(\ref{subsec:ContinuationMethods}). Bifurcation analysis reveals the
relation between the stability of a solution and the model parameters
and continuation methods calculates approximate solutions in function
of a model parameter.  In section \ref{sec:Results} we show
the results of our simulations. In section
\ref{sec:ConclusionAndDiscussion} we conclude and give an outlook.
%
%
\subsection{Description of the mathematical model}
\label{subsec:DescriptionOfTheMathematicalModel}
Before constructing a compartmental model that describes the
concentration of growth hormones per cell and its transport through a
plant organ such as for example a leaf, its geometry must be
specified.

\subsubsection{Geometry of the cells}
\label{subsubsec:GeometryOfTheCells}
\begin{figure}
	\centering
	\includegraphics[width=\textwidth]{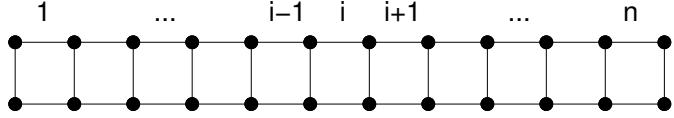}
	\caption{ \change{The one-dimensional domain, a row of regular square cells from $1$ to $n$.}         
          }
  \label{fig:domain} 
\end{figure}

\change{The domain in this work is a regular one dimensional row of
  cells, as in figure \ref{fig:domain}. Each cell is labeled and the
  set of cells is denoted with $V$. Therefore, for every cell in $V$,
  we can define the neighboring cells, a subset of $V$. For example,
  $\mathcal{N}_i = \left\{i-1, i+1\right\}$ is the set of neighboring
  cells of cell $i$. Further every cell consists of a number of cell
  walls.  The length of a cell wall between cell $i$ and cell $j$ is
  denoted with $l_{ij} = l_{ji}$.}

\change{The one dimensional domain in figure \ref{fig:domain}
  represents the 1D geometry of a group of cells in a part of a plant
  organ at a certain moment in time. In advanced models this geometry
  changes over time because cell walls grow and cells divide
  \citep{Smith2006}, but in this paper we look at the basic, coupled
  cell polarization and auxin transport model of Smith \textit{et al.}
  and the geometry is assumed to be static. As a consequence the
  row of cells is fixed.  The length of the cell wall is taken to be
  the length unit. As a result the volume of a cell also is the unit volume.}

\subsubsection{Transport of growth hormones}
\label{subsubsec:TransportOfGrowthHormones}
\change{For this} geometry it is now possible to formulate a model
for the transport of the growth hormones through the cells. We will
\change{write down with a concise mathematical notation the coupled
equations} of \citet{Smith2006} applied to a section across a leaf
with the geometry specified above. \change{Subsequently we introduce a 
slight generalization of the active transport term in the equation. }

In every cell $i$ two substances play an important role in the growth
process:
\begin{itemize}
\item The concentration of proteins PIN$1$ in cell $i$, which is time
  dependent, is denoted as \change{$\left[\textrm{PIN}1\right]_i\left(t\right) \in \mathbb{R}^+$ and is
  measured in micromol, i.e 10$^{-6}$mol, over the unit volume. We use the dimensionless variable $p_i\left(t\right) = \left[\textrm{PIN}1\right]_i\left(t\right)/\mu$M.}
\item The concentration of the hormone IAA in cell $i$, also known as
  auxin, is also time dependent and is denoted as \change{$ 
  \left[\textrm{IAA}\right]_i\left(t\right) \in \mathbb{R}^+$. Again in units of $\mu$M.
  We use the dimensionless variable $\left[\textrm{IAA}\right]_i\left(t\right)/\mu$M, which is denoted
  as $a_i\left(t\right)$. }
\end{itemize}
The model describes the evolution of these \change{dimensionless} concentrations in each
cell. This evolution depends in a non-linear way on the concentrations
of the neighboring cells. The value of $p_i\left(t\right)$ is
determined by the production and decay of PIN$1$. Its time evolution
for each cell $i \in V$ is modeled by
\begin{equation}
	\label{eq:PIN_p}
		\dfrac{dp_i\left(t\right)}{dt} = \dfrac{\rho_{_{\PIN_0}}+ \rho_{_{\PIN}}a_i\left(t\right)}{1 +  \kappa_{_{\PIN}}p_i\left(t\right)} - \mu_{_{\PIN}}p_i\left(t\right),
\end{equation}
where $\rho_{{_{\PIN}}_0} \in \mathbb{R}^+$ is the base production of
PIN1 proteins and $\rho_{_{\PIN}}\in \mathbb{R}^+$ is a coefficient
capturing the up-regulation of PIN1 production by auxin, \change{both
  measured per second}, $\kappa_{_{\PIN}}\in \mathbb{R}^+$
is the saturation coefficient of the PIN1 production, \change{which is
  dimensionless}, and $\mu_{_{\PIN}}\in \mathbb{R}^+$ is the PIN1
decay constant, \change{which has units of $1$/s.} This means that the
evolution of $p_i\left(t\right)$ in time depends strictly on the
concentration in the cell itself.

The concentration of IAA in a cell depends \change{not only on}
the production and the decay of auxin in the cell. Change of
$a_i\left(t\right)$ is also determined by diffusion (passive
transport) and active transport of auxin between the cells. The
change over time of the concentration of IAA is modeled by the
equation
\begin{eqnarray}
	\label{eq:IAA_a}
        \dfrac{da_i\left(t\right)}{dt} &=& \dfrac{\rho_{_{\IAA}}}{1+\kappa_{_{\IAA}}a_i\left(t\right)} - \mu_{_{\IAA}}a_i\left(t\right) - \sum_{j \in \mathcal{N}_i} \Dn\left(a_i\left(t\right) - a_j\left(t\right)\right) \\
        & & \left.\right. \hspace{0.5 cm}+ \sum_{j \in \mathcal{N}_i} \left(\text{ActiveTransport}_{j\rightarrow i} - \text{ActiveTransport}_{i \rightarrow j} \right), \notag
\end{eqnarray}
where $\rho_{_{\IAA}}\in \mathbb{R}^+$ is the IAA production
coefficient \change{which is measured per second},  
$\kappa_{_{\IAA}}\in \mathbb{R}^+$ is the dimensionless coefficient
which controls the saturation of IAA production, $\mu_{_{\IAA}}\in
\mathbb{R}^+$ is the IAA decay constant and $D\in \mathbb{R}^+$ is the
IAA diffusion coefficient, \change{both measured per second.} The active transport depends on the presence of
PIN1 denoted by $p_i$ and is modeled by the formula
\begin{equation}
\label{eq:ActiveTransport}
\text{ActiveTransport}_{i \rightarrow j} = \Tn \left(\dfrac{p_i\left(t\right)l_{ij} \textnormal{b}^{a_j\left(t\right)}}{\sum_{k\in \mathcal{N}_i}{l_{ik}\textnormal{b}^{a_k\left(t\right)}}}\right)	
\dfrac{a_i\left(t\right)^2}{1+\kappa_{_{\T}}a_j\left(t\right)^2},
\end{equation}
where $T\in \mathbb{R}^+$ is a polar IAA transport coefficient \change{expressed per second}, $b \in
\mathbb{R}^+$ is the exponentiation base which controls the extent to
which the PIN1 protein distribution is affected by the neighboring
cells and $\kappa_{_{T}} \in \mathbb{R}^+$ is an IAA transport
saturation coefficient. \change{These parameters are dimensionless.}
From equation \eqref{eq:IAA_a} we know that
the evolution of $a_i(t)$ depends only on itself, the first and the
second \change{nearest} neighbors of cell $i$. \change{Since we can specify the
  neighbors for every cell, the second nearest neighbors can be easily
  determined.  For example the second neighbors of cell $i$ are all
  elements in $\bigcup_{j \in \mathcal{N}_i} \mathcal{N}_j$.}
\change{Remark that because the length of each cell wall is taken to be the length unit, it cancels
from the equation.}

Equations \eqref{eq:PIN_p}, \eqref{eq:IAA_a} and
\eqref{eq:ActiveTransport} describe the basic coupled model of Smith
\textit{et al.}  that has been used to study the transport of hormones
in the Arabidopsis shoot apex.  \change{It differs mainly from
  other transport models by the active transport term. Smith
  \textit{et al.} uses a quadratic dependence to describe the flux on
  the auxin concentrations instead of a linear dependence. Further
  they introduce an exponential dependence of the localization of PIN1
  on the concentration of IAA. Therefor we will also consider two
  other models. One where the active transport is modeled with a
  linear dependence to describe the flux on the auxin concentration
\begin{equation}
	\label{eq:ActiveTransport2}
	\text{ActiveTransport}_{i \rightarrow j} = \Tn \left(\dfrac{p_i\left(t\right) \textnormal{b}^{a_j\left(t\right)}}{\sum_{k\in 																																									\mathcal{N}_i}{\textnormal{b}^{a_k\left(t\right)}}}\right)	
																								\dfrac{a_i\left(t\right)}{1+\kappa_{_{\T}}a_j\left(t\right)},
\end{equation}
and one model without the exponential dependence of the localization of PIN1 on the concentration of IAA
\begin{equation}
	\label{eq:ActiveTransport3}
	\text{ActiveTransport}_{i \rightarrow j} = \Tn \left(\dfrac{p_i\left(t\right) a_j\left(t\right)}{\sum_{k\in \mathcal{N}_i}{a_k\left(t\right)}}\right)	
																							\dfrac{a_i\left(t\right)^2}{1+\kappa_{_{\T}}a_j\left(t\right)^2}.
\end{equation}
The three different equations that model the active transport (equations \eqref{eq:ActiveTransport}, \eqref{eq:ActiveTransport2} and \eqref{eq:ActiveTransport3}) can be combined in one generalized equation
\begin{equation}
	\label{eq:ActiveTransport_all}
	\text{ActiveTransport}^{\omega,\tau}_{i \rightarrow j} = \Tn \left(\dfrac{p_i\left(t\right) \left(\omega\textnormal{b}^{a_j\left(t\right)} + \left(1-\omega\right)a_j\left(t\right)\right)}{\sum_{k\in 																																												\mathcal{N}_i}{\left(\omega\textnormal{b}^{a_k\left(t\right)} + 	\left(1-\omega\right)a^k\left(t\right)\right)}}\right)	
																													\dfrac{a_i\left(t\right)^\tau}{1+\kappa_{_{\T}}a_j\left(t\right)^\tau},
\end{equation}
For $\omega = 1$ and $\tau = 2$, \eqref{eq:PIN_p}, \eqref{eq:IAA_a} and
\eqref{eq:ActiveTransport_all} describe the equations of \citet{Smith2006}.}

\begin{figure}
	\centering
	\includegraphics[width=0.6\textwidth]{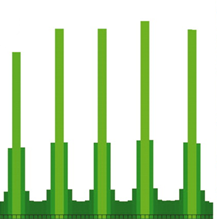}
	\caption{ \change{Simulation results of \citet{Smith2006} for a one
  dimensional row of $50$ cells with periodic boundary conditions. The shades of green show
  a difference in the IAA concentration. (figure reproduced from \citet{Smith2006}).}}
\label{fig:smithresults} 
\end{figure}

In \citet{Smith2006} a row of 50 equal
sized cells with periodic boundary conditions \change{was investigated}.  \change{The results of the} time evolution,
starting from an initially flat solution with a small amount of noise
to break symmetry, showed the emergence of a pattern in the IAA
concentrations (figure \ref{fig:smithresults}). Some cells have a 
very high concentration of \change{auxin}. \change{It was found that the peaks in a pattern are equally spaced and become} more prominent for an increasing IAA transport
coefficient T.

%
%
%
%
\section{The simulation problem}
\label{sec:TheSimulationProblem}

%
%
\subsection{Domain, boundary conditions and parameters}
\label{subsec:TheDomainItsBoundaryConditionsAndTheParameters}
In this paper we analyze the solutions of equations \eqref{eq:PIN_p}
and \eqref{eq:IAA_a}, for a
one dimensional file of equal sized square cells.  
\change{We assume that this file of cells represents a part of a leaf from the left margin to the midvein. This assumption is necessary in order to specify the boundary condition. Other parts of plant organs can give rise to other boundary conditions which will result in a small change in the model.}

\begin{figure}
	\centering
	\includegraphics[width=\textwidth]{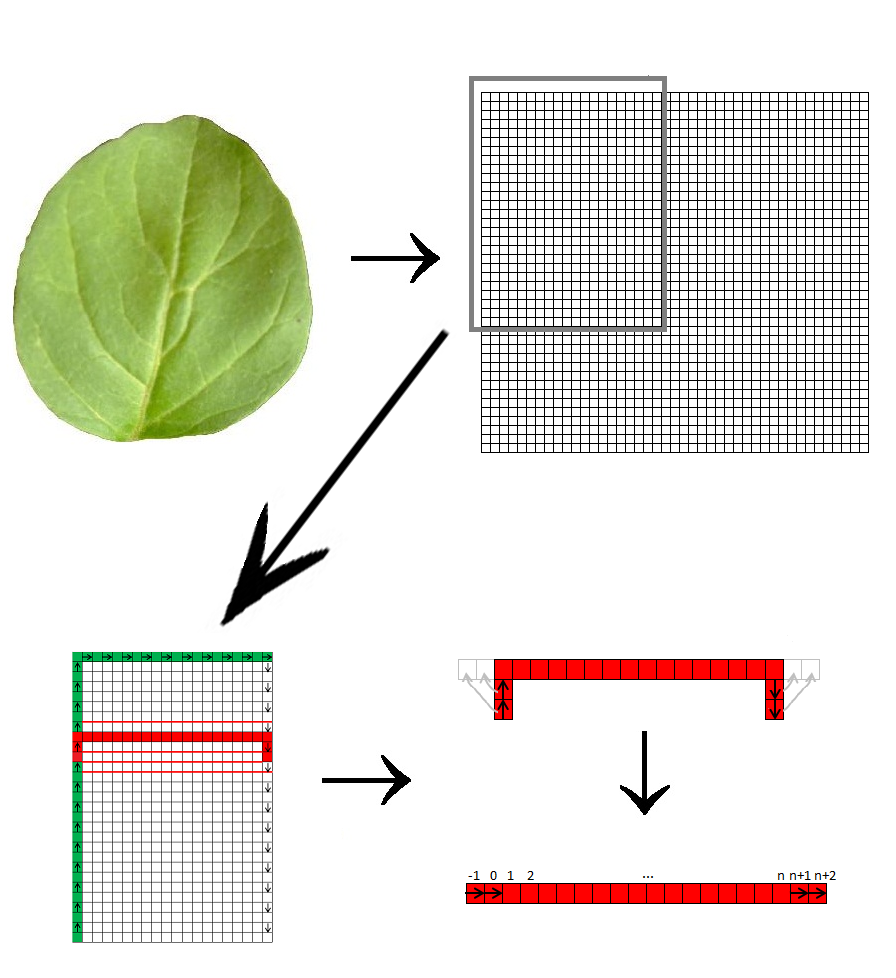}
	\caption{ 
					The one-dimensional
          model can be seen as a row of cells cut out of a leaf with
          equal sized square cells. In the first step we represent the
          leaf as a two dimensional squared grid of equal sized square
          cells. The second arrow indicates that we only consider the
          part from the left margin till the midvein. Here in each
          cell at the boundary, the direction of the auxin fluxes is
          indicated. In the third step a horizontal row of inner cells
          with, at each side two boundary cells is cut out of the
          domain. In the last step this domain is enrolled so it forms
          a one dimensional file of equal sized square cells. The two
          boundary cells at each side of the domain are the ghost
          cells.  
          }
  \label{fig:domainsteps} 
\end{figure}

To provide the boundary conditions, two ghost cells are required at
each end of the domain, since the model relates each cell with two
cells at the left and the right (figure \ref{fig:domainsteps}).

\change{The $n$ interior cells are labeled $1$ to $n$. The ghost cells
are cells $-1$ and $0$ at the left of the domain and cells $n+1$ 
and $n+2$ at the right.}
The concentration of the IAA hormone in the two ghost cells on each
side are chosen to describe the influx at the boundary of the leaf and
the efflux at the vein. The IAA concentration then changes linearly at
the boundaries as if Neumann boundary conditions are applied. 
We assume zero flux boundary conditions. This means that the boundary conditions
become
\change{
\begin{equation}
  \label{eq:equalNeumannbc}
  \begin{cases}
    a_{-1}\left(t\right)  =  a_1\left(t\right)  \quad \text{and} \quad a_{0}\left(t\right)  =  a_1\left(t\right), \\ 
    a_{n+1}\left(t\right)  =  a_{n}\left(t\right) \quad \text{and} \quad a_{n+2}\left(t\right) =  a_{n}\left(t\right) .
  \end{cases}
\end{equation}
} The value of $p_{0}\left(t\right)$ and $p_{n+1}\left(t\right)$ in
the ghost cells is determined by equation \eqref{eq:PIN_p} that
couples it to the value of $a_i\left(t\right)$ in the ghost cell.
Note that $p_{-1}\left(t\right)$ and $p_{n+2}\left(t\right)$ do not
appear in the problem since equation \eqref{eq:ActiveTransport} does
not require it. Together with an initial condition, the problem is
transformed in an initial value problem that we can solve numerically
with a time step method.  Remark that these homogeneous Neumann
boundary conditions are different from periodic boundary
conditions. The concentrations in the cells on the left side of the
domain can indeed be different from the concentrations in the cells on
the right side of the domain.

Equations \eqref{eq:PIN_p}, \eqref{eq:IAA_a} and
\eqref{eq:ActiveTransport} contain $11$ parameters. A short
description can be found in table \ref{table:parametervalues} and
further details can be found in \citet{Smith2006}. The values of these
parameters must be real and positive. For the simulations in this paper
we used three different parameter sets, M1, M2 and M3.  Parameter set
M2 corresponds with the values used by \citet{Smith2006}.  Parameter set M1 and M3 contain the same values
for the parameters as set M2 except for the IAA production
coefficient.  The value of this parameter is higher in parameter set
M1 and lower in set M3.
\begin{table}
\center
\caption{Values of the parameters of equations \eqref{eq:PIN_p}--\eqref{eq:ActiveTransport} used in the simulations.  The parameter values of M2 are found in \citet{Smith2006}.}
\label{table:parametervalues}
\begin{tabular}[H]{|c|c|c|c|c|}
\hline
Symbol & Description & \multicolumn{3}{|c|}{Value} \\
\cline{3-5} 
 			 &						 & M1 & M2 &M3 \\
\hline \hline
b & Base for exponential PIN allocation & $3.000$ & $3.000$ & $3.000$\\
\hline
$\kappa_{_{\PIN}}$ & PIN saturation coefficient & $1.000$ & $1.000$ & $1.000$\\
\hline
$\kappa_{_{\T}}$ & Transport saturation coefficient & $1.000$ & $1.000$ & $1.000$\\
\hline
$\kappa_{_{\IAA}}$ & IAA saturation coefficient & $1.000$  & $1.000$ & $1.000$\\
\hline
$\rho_{_{\PIN_0}}$ & Base production of PIN & $0.000$ & $0.000$ & $0.000$\\
\hline
$\rho_{_{\PIN}}$ & PIN production coefficient & $1.000$ & $1.000$ & $1.000$ \\
\hline
$\mu_{_{\PIN}}$ & PIN decay coefficient & $0.100$ & $0.100$ & $0.100$\\
\hline
$\mu_{_{\IAA}}$ & IAA decay coefficient & $0.100$ & $0.100$ & $0.100$\\
\hline
$\rho_{_{\IAA}}$ & IAA production coefficient & $1.500$ & $0.750$ & $0.500$\\
\hline
D & IAA diffusion coefficient & $1.000$ & $1.000$ & $1.000$\\
\hline
T & IAA transport coefficient & $3.500$  & $3.500$  & $3.500$ \\
\hline
\end{tabular}
\end{table}

%
%
\subsection{Time integration}
\label{subsec:TimeIntegration} 
Similar to \citet{Smith2006}, we can solve the initial value
problem of Smith and collaborators with numerical integration. Analysis of the eigenvalues of the
Jacobian shows that they are mostly located along the negative real
axis with some small complex conjugate pairs of outliers at the left of the
imaginary axis. This suggests that the fourth order Runge-Kutta method
\citep{Hairer2009} with time step $\Delta t = 0.01$ results in a stable
method to integrate the equations.
\begin{figure}
  \centering
	\includegraphics[width=\textwidth]{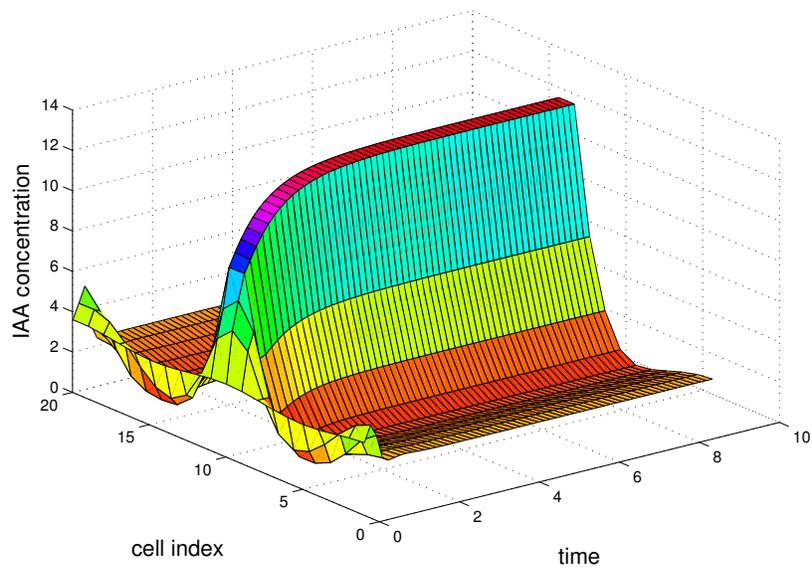}
        \caption{The time evolution of the solution \change{of equations \eqref{eq:PIN_p}, \eqref{eq:IAA_a} and \eqref{eq:ActiveTransport}} for a row of 20
          cells with zero Neumann boundary condition and parameter set
          M1. The initial condition is given by equation
          \eqref{eq:beginvalue1} and we used RK4 for numerical
          integration.}
	\label{fig:timeEvolution}
\end{figure}
In figure \ref{fig:timeEvolution} the time evolution is shown from
$t=0$ to $t=10$. The domain contains $20$ cells plus $4$ ghost cells
where we assume zero Neumann boundary conditions. The parameter
values of set M1 are used and the initial value for the concentration is
\begin{equation}
\label{eq:beginvalue1}
	p_i(t=0) = 5.4  \quad \text{and} \quad   a_i(t=0) =3.4,
\end{equation}
\normalsize where a small perturbation
$0.2\sin\left(\left(5\left(i+2\right)\pi\right)/24\right)$ for
$i=1,\ldots, 20$
was added to it to break symmetry. \change{Any other initial state nearby will lead to the same long term solution.} \\
Figure \ref{fig:timeEvolution} shows the development of a pattern in
the concentration of \change{auxin}. After a certain time the pattern
arrives in a stable steady state. For the row of $20$ cells a single
peak with a high IAA concentration is formed.

%
%
\subsection{Steady state problem}
\label{subsec:SteadyStateProblem}
Rather than evolving the system in time, we can calculate the steady
state solutions directly. 
We rewrite equations \eqref{eq:PIN_p},
\eqref{eq:IAA_a} and \eqref{eq:ActiveTransport} in order to obtain the
steady state equations for this specific geometry and boundary
conditions. 
The steady state problem
becomes:
\begin{equation}
\label{eq:Neumannsteadystate}
\begin{cases}
  0 &= \dfrac{\rho_{_{\PIN_0}}+ \rho_{_{\PIN}}a_i}{1 +  \kappa_{_{\PIN}}p_i} - \mu_{_{\PIN}}p_i \quad \text{with} \quad i=0 ,1, \ldots , n,n+1 \\
  0 &= \dfrac{\rho_{_{\IAA}}}{{1 + \kappa_{_{\IAA}}}a_i} - \mu_{_{\IAA}}a_i	- {\Dn\left(a_i-a_{i-1}\right)} - {\Dn\left(a_i-a_{i+1}\right)} \\
  & +\Tn \left(\dfrac{p_{i-1}\textnormal{b}^{a_i}}{\textnormal{b}^{a_{i-2}} + \textnormal{b}^{a_{i}}}\right) \dfrac{a_{i-1}^2}{1+\kappa_{_{T}}a_i^2} - \Tn  \left(\dfrac{p_{i}\textnormal{b}^{a_{i-1}}}{\textnormal{b}^{a_{i-1}} +    \textnormal{b}^{a_{i+1}}}\right)
  \dfrac{a_{i}^2}{1+\kappa_{_{T}}a_{i-1}^2}\\
  & + \Tn \left(\dfrac{p_{i+1}\textnormal{b}^{a_i}}{\textnormal{b}^{a_{i}} + \textnormal{b}^{a_{i+2}}}\right) \dfrac{a_{i+1}^2}{1+\kappa_{_{T}}a_i^2}  - \Tn \left(\dfrac{p_{i}\textnormal{b}^{a_{i+1}}}{\textnormal{b}^{a_{i-1}} + \textnormal{b}^{a_{i+1}}}\right) \dfrac{a_{i}^2}{1+\kappa_{_{T}}a_{i+1}^2}   \\
  & \hspace{7cm} \textrm{with } \quad i = 1, \ldots , n\\ \\
 & a_{-1} = a_1  \quad \text{and} \quad  a_{0} = a_1  \\
 & a_{n+1} = a_{n} \quad \text{and} \quad  a_{n+2} = a_{n} ,
\end{cases}
\end{equation}
where the indices $0$ and $n+1$ in the first equation express the
  coupling of $p_i$ to $a_i$ in the first ghost cells.  This system can be written
as the system of equations
\begin{equation}
\label{eq:ContEq}
F\left(\mathbf{U}, \boldsymbol{\lambda}\right) = \mathbf{0} ,
\end{equation}
where $F: \R^{2n+m} \longrightarrow \R^{2n} : \left(\mathbf{U},
  \boldsymbol{\lambda} \right) \mapsto
F\left(\mathbf{U},\boldsymbol{\lambda} \right)$ with $n$ the number of
cells and $m$ the number of parameters. $\mathbf{U}$ is a $2n$
dimensional solution vector of the problem that contains both the $p$
and $a$ steady state variables and $\boldsymbol{\lambda} \in \R^m$ denotes the set
of parameters.

\change{For the equations \eqref{eq:PIN_p} and \eqref{eq:IAA_a} with active transport equations \eqref{eq:ActiveTransport2}, \eqref{eq:ActiveTransport3} or \eqref{eq:ActiveTransport_all} the steady state problem can be obtained in the same way.}

%
%
%
%
\section{The trivial solution}
\label{sec:TheTrivialSolution}
In this section we search for a trivial solution of system
\eqref{eq:Neumannsteadystate}, a solution that can be calculated
analytically and that will be used as a starting point for the numerical
continuation in section \ref{sec:Results}.

\begin{figure}
	\center
	\includegraphics[width = 0.9\textwidth]{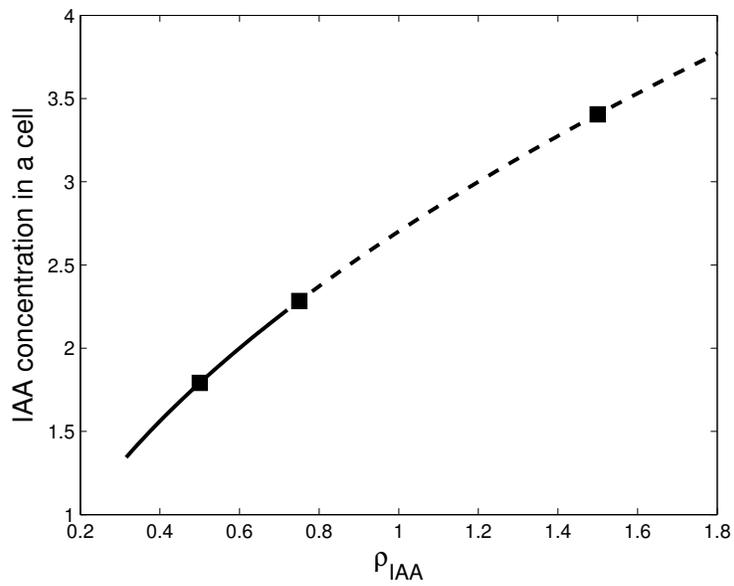} 
	\caption{ The concentration of IAA in the trivial
          solution as a function of the IAA production coefficient $\rho_{\IAA}$. Each cell has
          the same concentration. However, for large $\rho_{\IAA}$
          this solution becomes unstable. Stability is calculated \change{for equation \eqref{eq:Neumannsteadystate}} and a row of 20 cells. The other parameters are taken from table \ref{table:parametervalues}. }
	\label{fig:theorSol_rhoIAA}
\end{figure}

If we assume that the solution is homogeneous then the values of $p$
and $a$ are the same for all cells so that
\begin{equation} \label{eq:homogeneousConditionUV}
p_i = p_j  \quad \text{and} \quad  a_i = a_j   \qquad  \forall i,j =-1,\ldots,n+2. 
\end{equation}
The system \eqref{eq:Neumannsteadystate} now reduces to
\begin{equation}
	\begin{cases}
    0 = \dfrac{\rho_{_{\PIN_0}}+ \rho_{_{\PIN}}a_i}{1 +  \kappa_{_{\PIN}}p_i} - \mu_{_{\PIN}}p_i  \quad \text{ for } i=0,\ldots,n+1 \\
  	0 = \dfrac{\rho_{_{\IAA}}}{{1 + \kappa_{_{\IAA}}}a_i} - \mu_{_{\IAA}}a_i \quad \text{ for } i=1,\ldots,n ,
  	\label{eq:trivialSolutionequation}
  \end{cases}
\end{equation}
and
\begin{equation}
  a_{-1} = a_{0} = a_1 \quad \text{and} \quad a_{n+1} = a_{n+2} =  a_{n}.
\end{equation}
Because $p_i$ and $a_i$ are real positive numbers, \change{we find a unique} solution  \change{that is given by} 
\begin{equation}
  \begin{cases}
  p_i = \dfrac{-1 + \sqrt{1 + 4\kappa_{_{\PIN}}\left(\rho_{_{\PIN_0}} + \rho_{_{\PIN}}a_i\right)/\mu_{_{\PIN}}}}{2\kappa_{_{\PIN}}},  \\
  a_i = \dfrac{-1 + \sqrt{1 +
      4\kappa_{_{\IAA}}\rho_{_{\IAA}}/\mu_{_{\IAA}}}}{2\kappa_{_{\IAA}}}
  ,
  \label{eq:trivialSolution}
\end{cases}
\end{equation}
with $i=-1, \ldots, n+2$. This is the trivial solution of the system.
\change{Note that the same trivial solution is obtained for the models
with active transport equation   \eqref{eq:ActiveTransport2}, 
\eqref{eq:ActiveTransport3} or \eqref{eq:ActiveTransport_all} .}  
From equation
\eqref{eq:trivialSolutionequation} we know that for a certain
parameter set, there is only one trivial homogeneous solution.  By
formula \eqref{eq:trivialSolution} it is easy to calculate this
trivial solution for different parameter values.  Figure
\ref{fig:theorSol_rhoIAA} shows the concentration of IAA in one cell
(cell number 6) versus the parameter $\rho_{_{\IAA}}$. Because the
solution is homogeneous, the trivial solution curve would be the same
for every cell and is independent of the number of cells. Figure
\ref{fig:theorSol_rhoIAA} denotes also the trivial solution for
parameter set M1, M2 and M3 \change{with a square}.

%
%
\subsection{Stability of the trivial solution} 
\label{subsec:StabilityOfTheTrivialSolution} 
\change{The
  value of the IAA concentration in the trivial solution is independent of
  the IAA diffusion coefficient D and the IAA transport coefficient T. Furthermore, it
  does not depend on the type of model for the active transport. Any
  choice of $\omega$ or $\tau$ yields to the same trivial
  solution. However, the stability depends in a sensitive way on the
  diffusion and transport coefficients and on the model of the active
  transport.}  Although expression \eqref{eq:trivialSolution} for the
trivial solution is readily obtained, it is not easy to determine the
stability of this solution. The Jacobian matrix of the coupled system
\eqref{eq:Neumannsteadystate} is not trivial. \change{It has a sparse
  blocked structure where the first top-left block is diagonal and
  contains the linearization of Eq.~\eqref{eq:PIN_p}. The top-right
  block is also diagonal and has the derivative of
  Eq.~\eqref{eq:PIN_p} to $a_i$.  The bottom-left block is tridiagonal
  since the active transport in cell $i$ depends on $p_{i-1}$, $p_{i}$ and
  $p_{i+1}$.  Finally, the bottom right block has five diagonals because 
  the active transport depends on the neighbors and the neighbors of the neighbors.
  Because of the complicated structure the eigenvalues are non-trivial
  and we have found that for the coupled problem they are not
  necessarily real valued in contrast to the results of the uncoupled
  problem, see for example \citep{Jonsson2006}.}

The stability can, however, \change{easily} be calculated by numerical
means and in our simulations we approximate the Jacobian with central
finite differences. The $j$-th column of
$J(\mathbf{U},\boldsymbol{\lambda})$ is
$J(\mathbf{U},\boldsymbol{\lambda})\mathbf{e}_j$ where $\mathbf{e}_j$
is the unit vector with the $j$-th component \change{equal to 1} and
the other components \change{equal to} 0. The column is then
approximated as
\begin{equation}
J(\mathbf{U},\boldsymbol{\lambda})\mathbf{e}_j =  \frac{F(\mathbf{U} + \epsilon \mathbf{e}_j, \boldsymbol{\lambda}) - F(\mathbf{U} - \epsilon \mathbf{e}_j, \boldsymbol{\lambda})}{2\epsilon},
\end{equation}
where $\epsilon$ is taken of the order of 10$^{-7}$.  \change{Once an
approximation to the Jacobian is obtained, its eigenvalues can be
calculated.}  This numerical approach can also be used to study the
stability of other solutions.

The stability of the trivial solution  of equation \eqref{eq:Neumannsteadystate} 
for a row of $20$ cells is shown
in figure~\ref{fig:theorSol_rhoIAA}. For smaller values of
$\rho_{_{\IAA}}$ the eigenvalues of the trivial solutions lie in the
left half-plane of the complex plane. Therefor these \change{stable} solutions 
\change{are drawn} with a full line in figure
\ref{fig:theorSol_rhoIAA}. For larger values of $\rho_{_{\IAA}}$, at
least one eigenvalue lies in the right half plane and so the trivial
solution is unstable. This is indicated with a dotted line.

\begin{figure}
	\centering
	\subfloat[] {\label{fig:stable_unstable_region_rhoIAA_T} 	\includegraphics[width=0.33\textwidth]{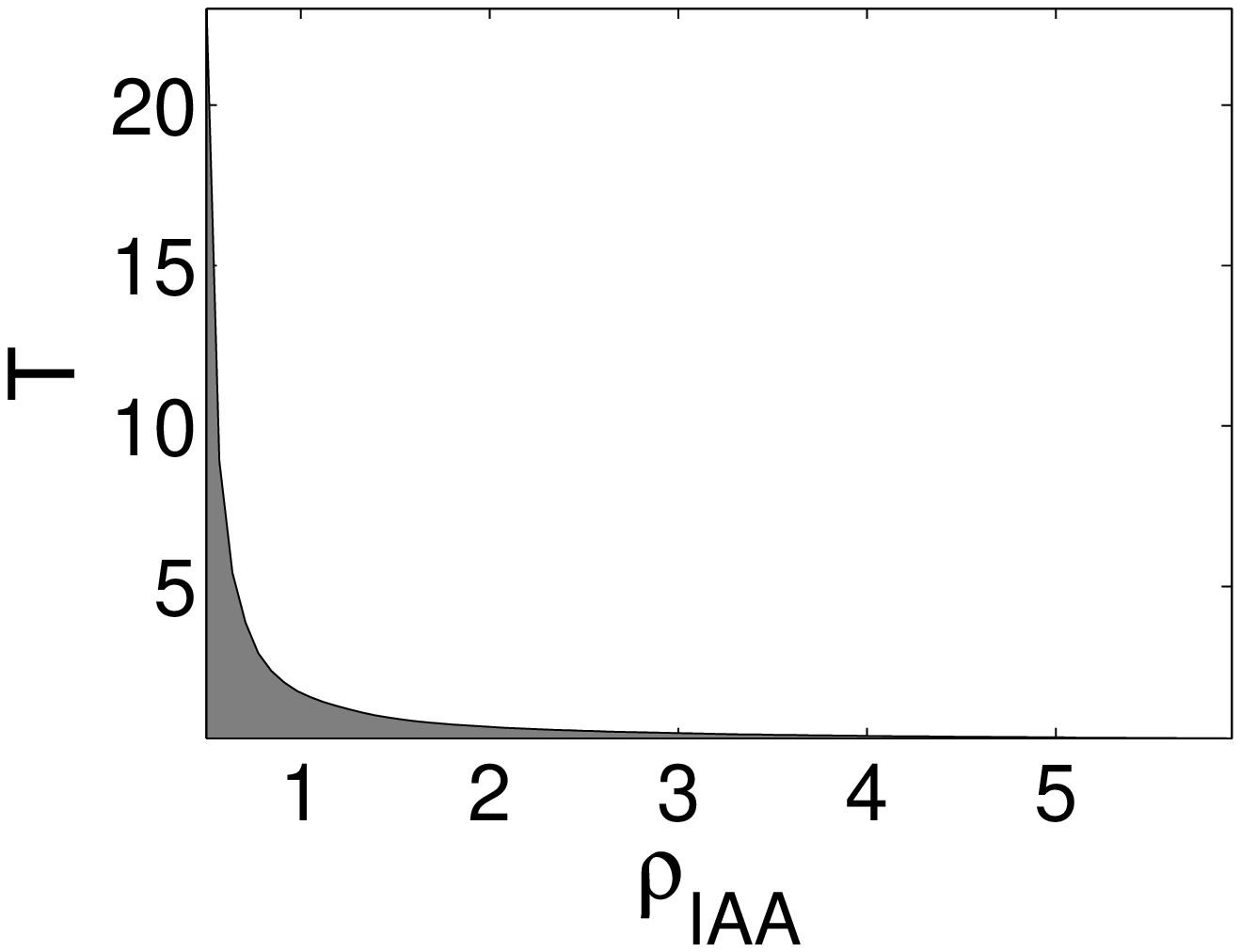}}
	\subfloat[] {\label{fig:stable_unstable_region_rhoIAA_D} 	\includegraphics[width=0.33\textwidth]{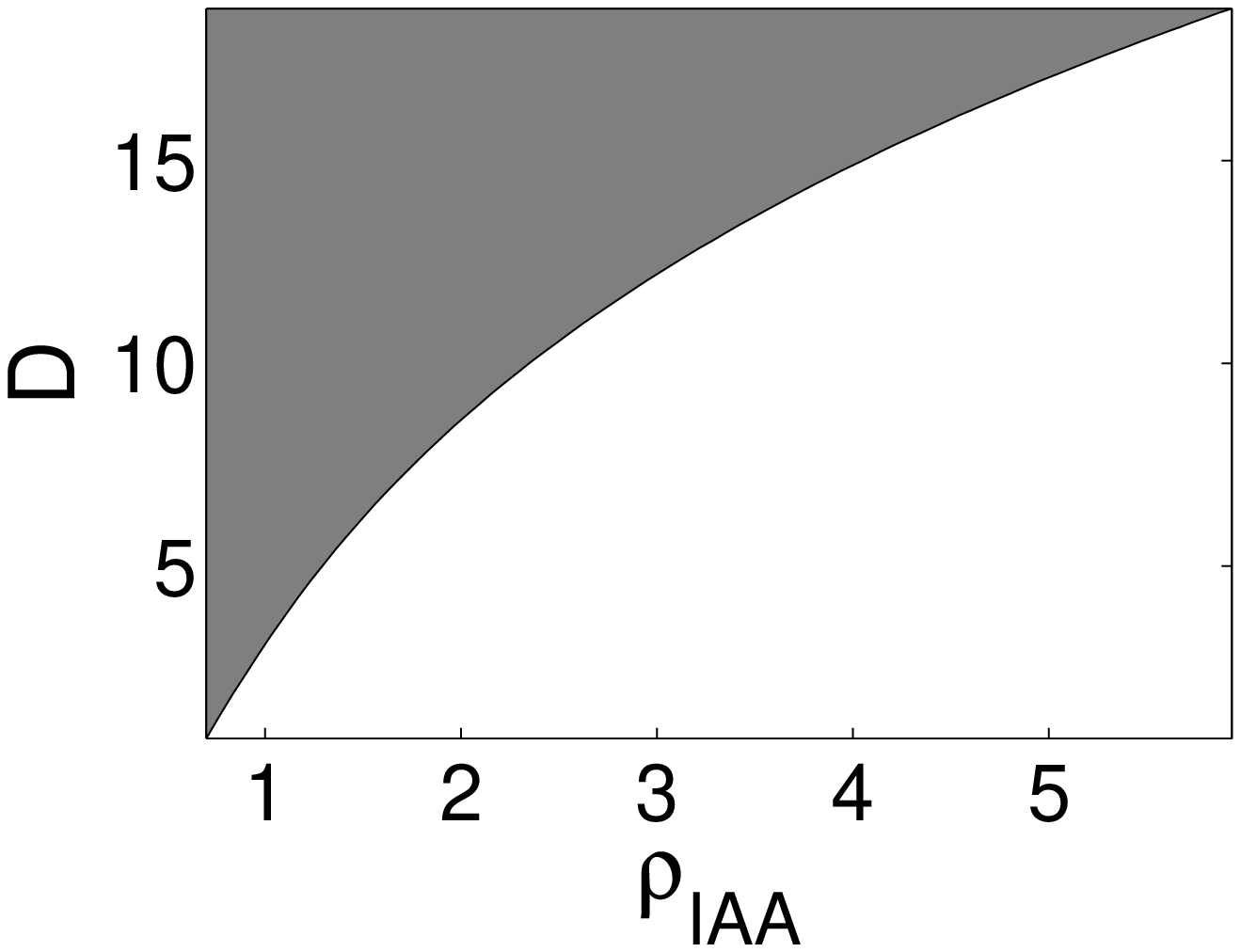}}
	\subfloat[] {\label{fig:stable_unstable_region_D_T}				\includegraphics[width=0.33\textwidth]{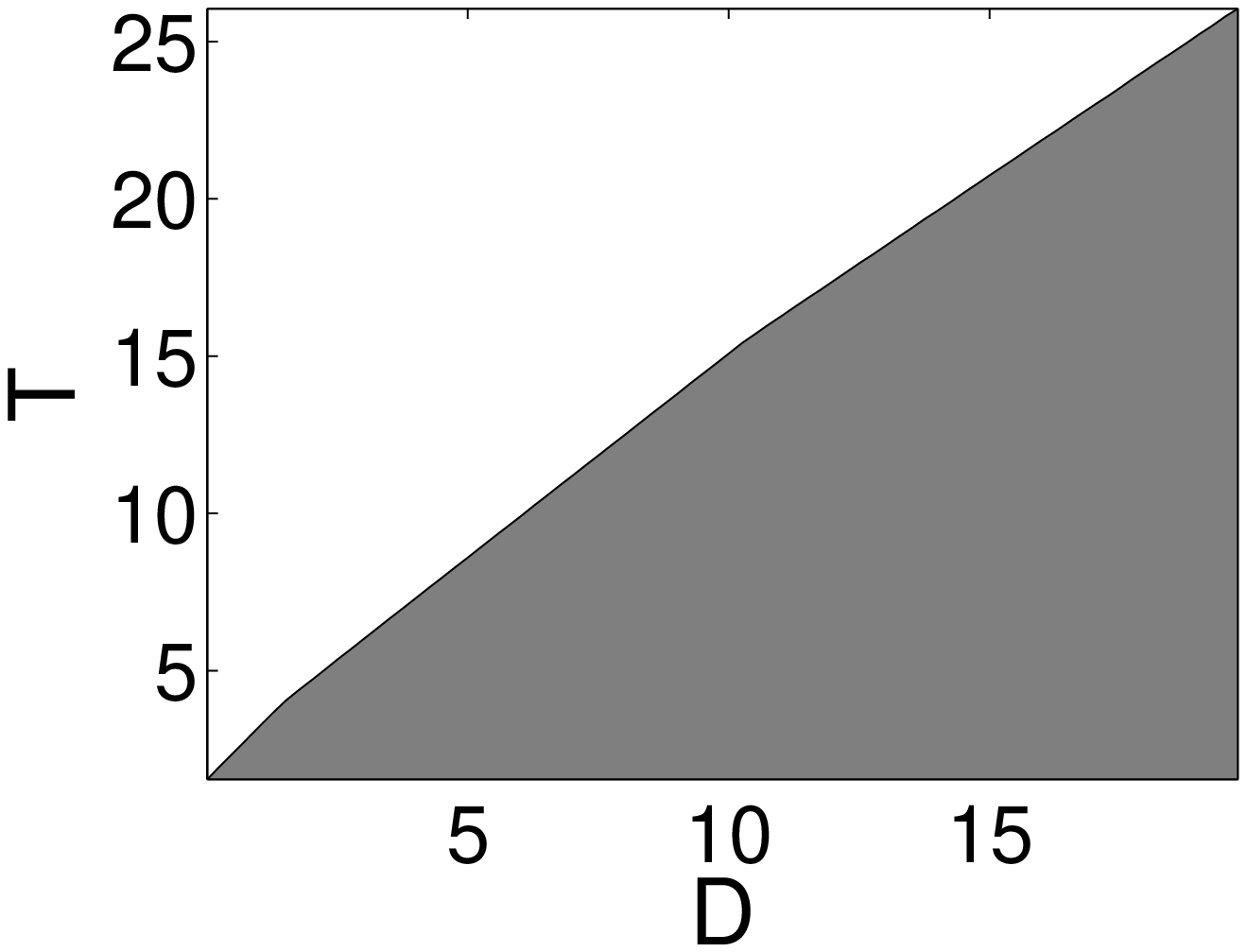}} \\
	\subfloat[] {\label{fig:stable_unstable_region_rhoIAA_T2} 	\includegraphics[width=0.33\textwidth]{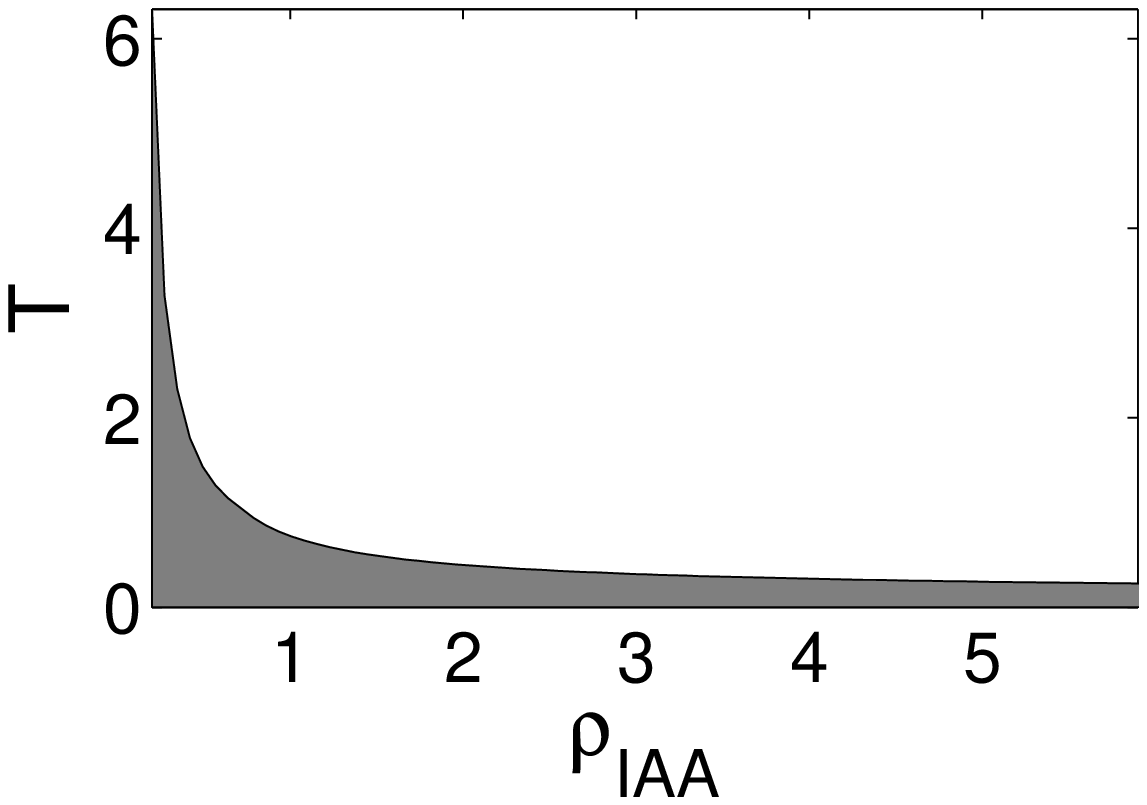}}
	\subfloat[] {\label{fig:stable_unstable_region_rhoIAA_D2} 	\includegraphics[width=0.33\textwidth]{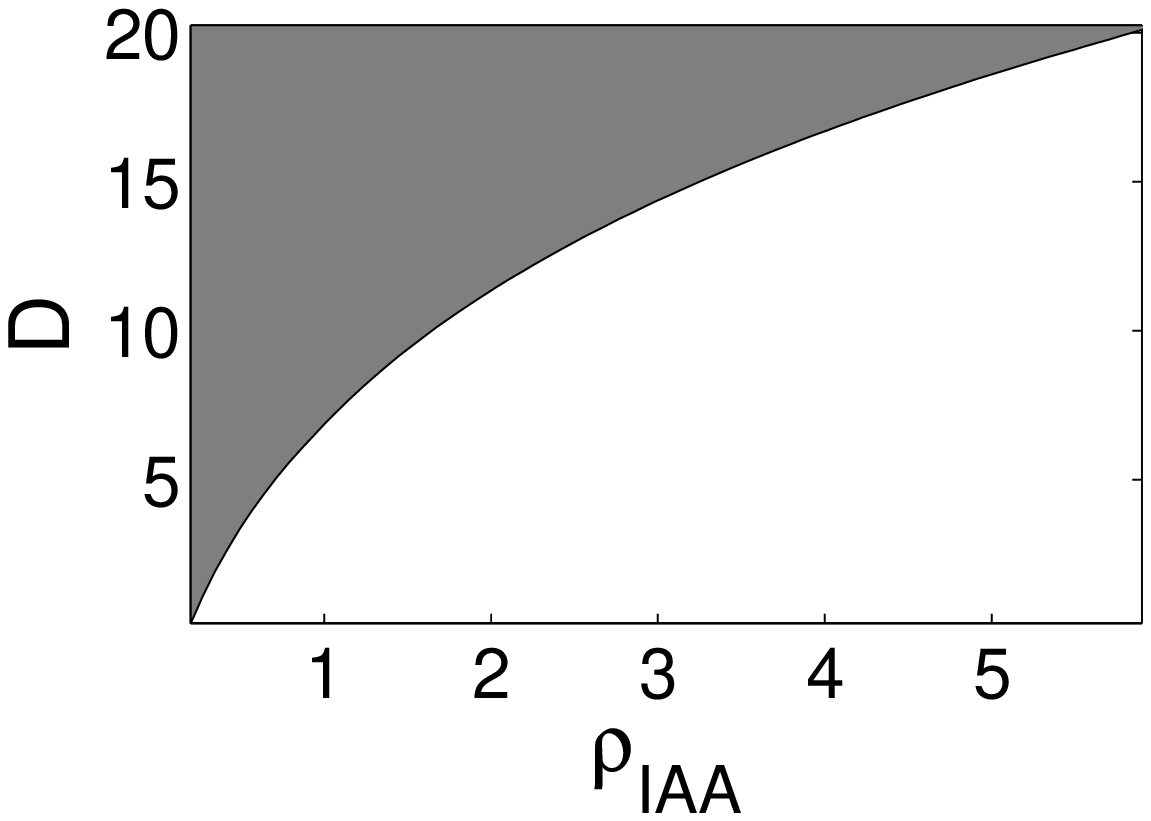}}
	\subfloat[] {\label{fig:stable_unstable_region_D_T2}				\includegraphics[width=0.33\textwidth]{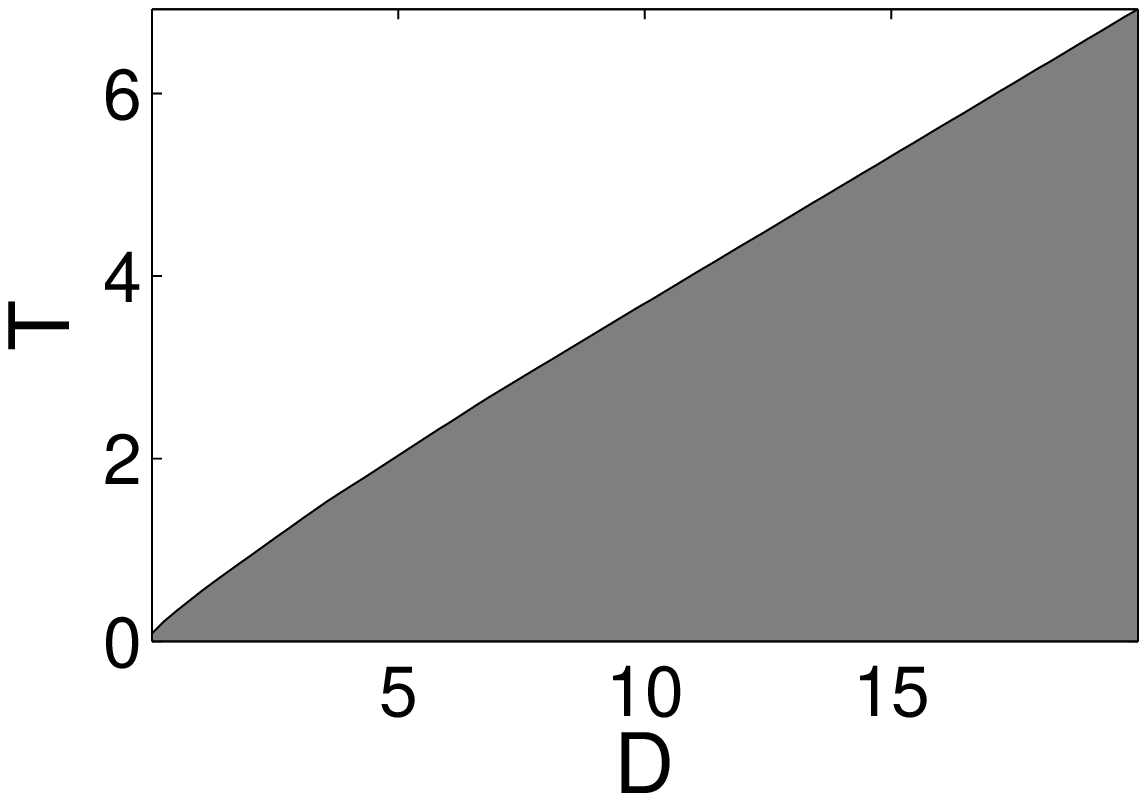}}
	
	\caption{Stable (marked in gray) and unstable region of the
          trivial solution for a row of 20 cells and different choices
          of the parameters $\rho_{\IAA}$ (IAA production
          coefficient), D (IAA diffusion coefficient) and T (IAA
          transport coefficient). Other parameters are taken from
          parameter set M2. \change{The first row of figures displays the 
          stability regions for the model with active transport equation \eqref{eq:ActiveTransport},
          the second row for the model with active transport equation \eqref{eq:ActiveTransport2}.}}
	\label{fig:stable_unstable_region}
\end{figure}

Also for other parameter values we can calculate the stability of the
trivial solution. In each plot on figure
\ref{fig:stable_unstable_region} two parameter values are varied. The
other parameter values are taken as in parameter set M2. 
\change{The first row of figures corresponds with the stability 
region for equation \eqref{eq:Neumannsteadystate} 
( i.e. the basic coupled model of \citet{Smith2006}) 
and the second row of figures represents the stability regions 
for the steady state problem with a modified active transport equation \eqref{eq:ActiveTransport2}.} 
These plots show where the trivial solution is stable (marked in gray) for a row of $20$ cells. 
For example, figure \ref{fig:stable_unstable_region_rhoIAA_T} shows
that very small values of $\rho_{_{\IAA}}$ give a stable trivial
solution for almost every value of T \change{in the original model.}
All other values of $\rho_{_{\IAA}}$ give an unstable trivial solution
if T is not too small.  \change{Further we find that when the active
  transport is modeled with a linear dependence on the auxin
  concentration, the shape of the stability region of the trivial
  solution remains approximately the same but it is much smaller than
  with a quadratic dependency (compare for example figure
  \ref{fig:stable_unstable_region_rhoIAA_T} and figure
  \ref{fig:stable_unstable_region_rhoIAA_T2}). For both models we find
  that increasing the number of cells, leads to approximately the same
  shape of the stable region of the trivial solution -- it only gets
  slightly smaller.  Applying the model without the exponential
  dependence of the localization of PIN1 on the concentration of IAA
  results in a stable trivial solution for the entire tested range of
  parameters.}

%
%
%
%
\section{Methods}
\label{sec:Methods}

%
%
\subsection{Bifurcation analysis}
\label{subsec:BifurcationAnalysis}
The study of the relation between the stability of a solution and the
parameters of the corresponding dynamical system is known as
bifurcation analysis \citep{Seydel}.  Such an analysis identifies the
stable and unstable solutions and the bifurcation points that mark the
transitions between them.  This is biologically relevant since it will
allow us to predict the patterns that emerge in the time evolution as
the parameters of the model are changed.  A bifurcation point is a
solution $\left(\textbf{U}_i,\boldsymbol{\lambda}_i\right)$ of system
\eqref{eq:ContEq} where the number of solutions changes when
$\boldsymbol{\lambda}$ passes $\boldsymbol{\lambda}_i$.  In this
article there are several types of bifurcation points such branch
points, limit points and Hopf bifurcation points that will play a
role.  A branch point is a bifurcation point where two or more
branches with distinct tangents intersect. A limit point, also called
a turning point, is a point where, locally, no solutions exist on one
side of the limit point and two solutions on the other side. A Hopf
bifurcation is a transition where a periodic orbit appears and branch
points and limit points are both bifurcation points among steady state
solutions. For complete review of their properties wer refer to
\citep{Seydel}.  The analysis usually leads to a bifurcation diagram
that highlights the connections between stable and unstable branches
as the parameters change. It is useful to track all these solution
branches that emerge, split or end in a bifurcation point. This can be
done with the help of numerical continuation methods.

%
%
\subsection{Continuation methods}
\label{subsec:ContinuationMethods}
The system of equations \eqref{eq:ContEq} is a smooth map and we know
that $\mathbf{0} \in \text{Range}\left(F\right)$. Following the implicit
function theorem we know that for a regular point
$\mathbf{x}_0=\left(\mathbf{U}_0,\boldsymbol{\lambda}_0\right) \in
\mathbb{R}^{2n+m}$ of $F$ that satisfies $F\left(\mathbf{x}_0\right) =
\mathbf{0} $, the solution set $F^{-1}\left(\mathbf{0}\right)$ can be locally
parametrized about $\mathbf{x}_0$ with respect to some parameter
$s$. This means that the system of equations $F\left(\mathbf{U},
  \boldsymbol{\lambda}(s)\right) = \mathbf{0}$ defines an implicit curve
$\mathbf{U}(\boldsymbol{\lambda}(s))$ where
$\boldsymbol{\lambda}(s):\mathbb{R}\rightarrow \mathbb{R}^m$ is any
parametric curve in the $\mathbb{R}^m$ \citep{Allgower}. The idea of
continuation methods is to find a curve $c$ of approximate solutions
$\textbf{U}$ of the system in function of the parameter
$\boldsymbol{\lambda}(s)$. To construct such a curve of subsequent
solution points $\mathbf{x}_i =
\left(\mathbf{U}_i,\boldsymbol{\lambda}_i\right) =
\left(\mathbf{U}_i,\boldsymbol{\lambda}(s_i)\right)$, continuation
methods use a starting point $\mathbf{x}_0 =
\left(\mathbf{U}_0,\boldsymbol{\lambda}_0\right)$, a solution of
system \eqref{eq:ContEq}, along with an initial continuation
direction \citep{Krauskopf}. This starting point is typically a trivial
solution. An important family of the continuation methods are the
predictor-corrector methods such as pseudo-arc-length continuation. The
idea of the algorithm is to first predict a new solution point. In the
corrector step, this predicted point is the start value for an
iterative method that will approximate the solution to a given
tolerance. For the pseudo-arc-length, the predictor step uses the
tangent vector to the curve at a solution point and a given step size
to predict a guess for the next solution point on the curve. The
corrector step improves the guess with Newton iterations.

Numerical continuation is available in AUTO \citep{Doedel1997}, LOCA
part of Trilinos \citep{Salinger2005}, PyDS \citep{PyDS} and others.
These libraries can often also identify the bifurcations that occur
along the continued curve and some of them, such as AUTO can
automatically switch between branches at bifurcation points.

%
%
%
%
\section{Results}
\label{sec:Results}
\change{This section presents several examples that highlight specific
  properties of the dynamics of the model.  In the first three
  examples we give, for a file of $20$ cells, the numerical
  bifurcation analysis of equation \eqref{eq:Neumannsteadystate} with
  respectively parameter sets M1, M2 and M3. In the fourth example, 
  we enlarge the system and study now a
  file of $100$ cells with parameter set M1 instead of $20$ cells.
  In the last two examples we investigate the model with a generalized equation
  for the active transport. In example 5, we look at the influence of
  the parameter $\tau$ on the bifurcation scheme found in example
  $1$. Finally in example $6$ we investigate the influence of the
  parameter $\omega$ on the stability. }

\change{In examples $1$ to $5$ the IAA transport coefficient T is the continuation parameter.
 We have chosen T as the continuation parameter similar to the
one-dimensional simulations of \citet{Smith2006}. Also \citet{Jonsson2006}
investigated the influence of the IAA transport 
coefficient T in their simple model by changing the ratio D/TP, 
with P the fixed value for PIN1. In example 6 we use $\omega$ as continuation parameter.}

\change{Each time we choose to display a bifurcation diagram that depicts the 
IAA concentration in cell number $6$ versus the continuation 
parameter. Alternative choices for the measure on the y-axis (e.g. a different cell) would be equally valid. }

We will find that the trivial solution loses its stability through
either a branch point or a Hopf bifurcation.  
\change{The results also show the small effect of the quadratic 
dependence in comparison with the linear dependence to describe 
the flux on the auxin concentration.}

\begin{figure}
	\centering
	\subfloat[Bifurcation Diagram] {\label{fig:BifurcationDiagramScenario1} \includegraphics[width = 0.95\textwidth]{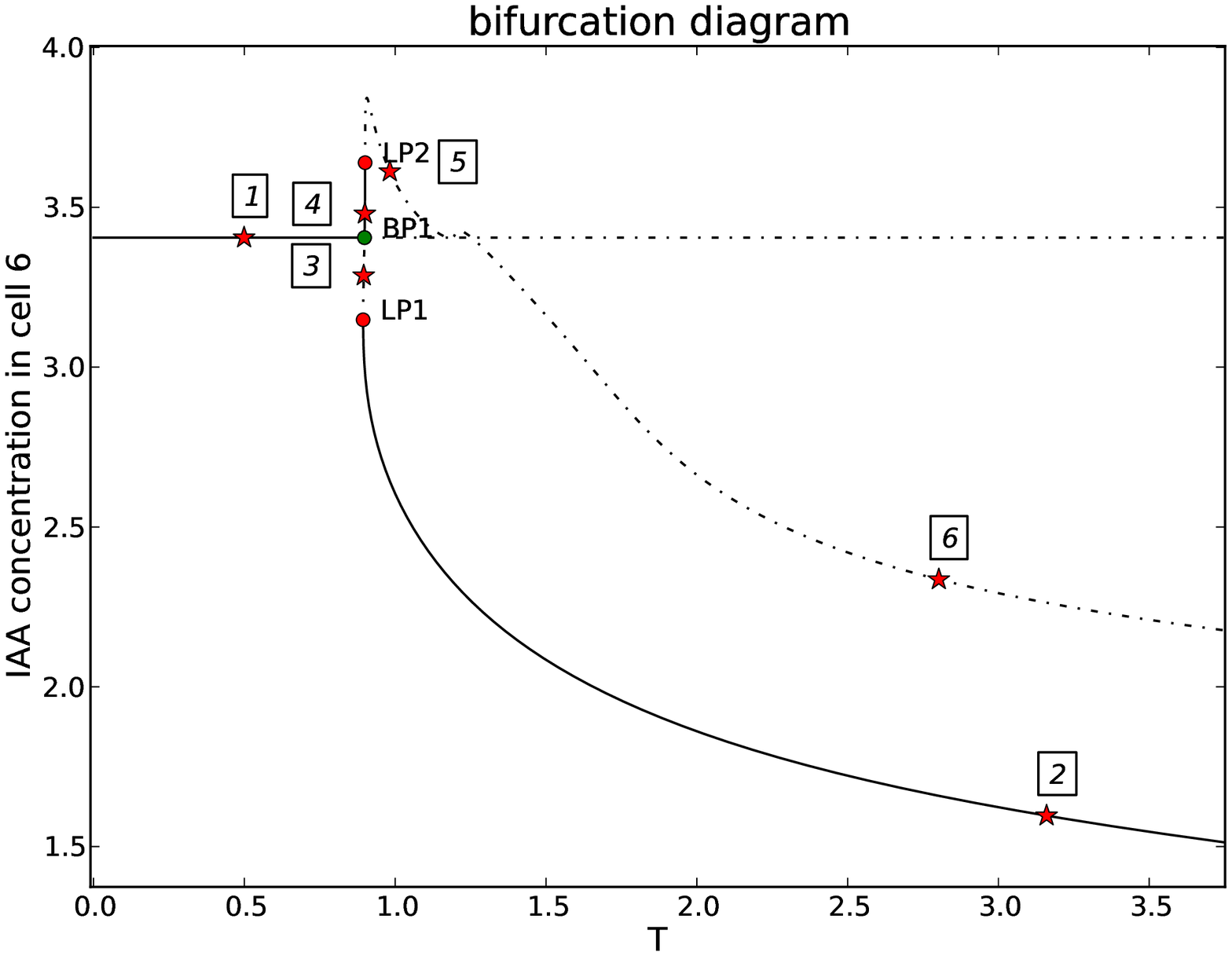}}\\
	\subfloat[Solution 1] {\label{fig:solution1}\includegraphics[width = 0.3\textwidth]{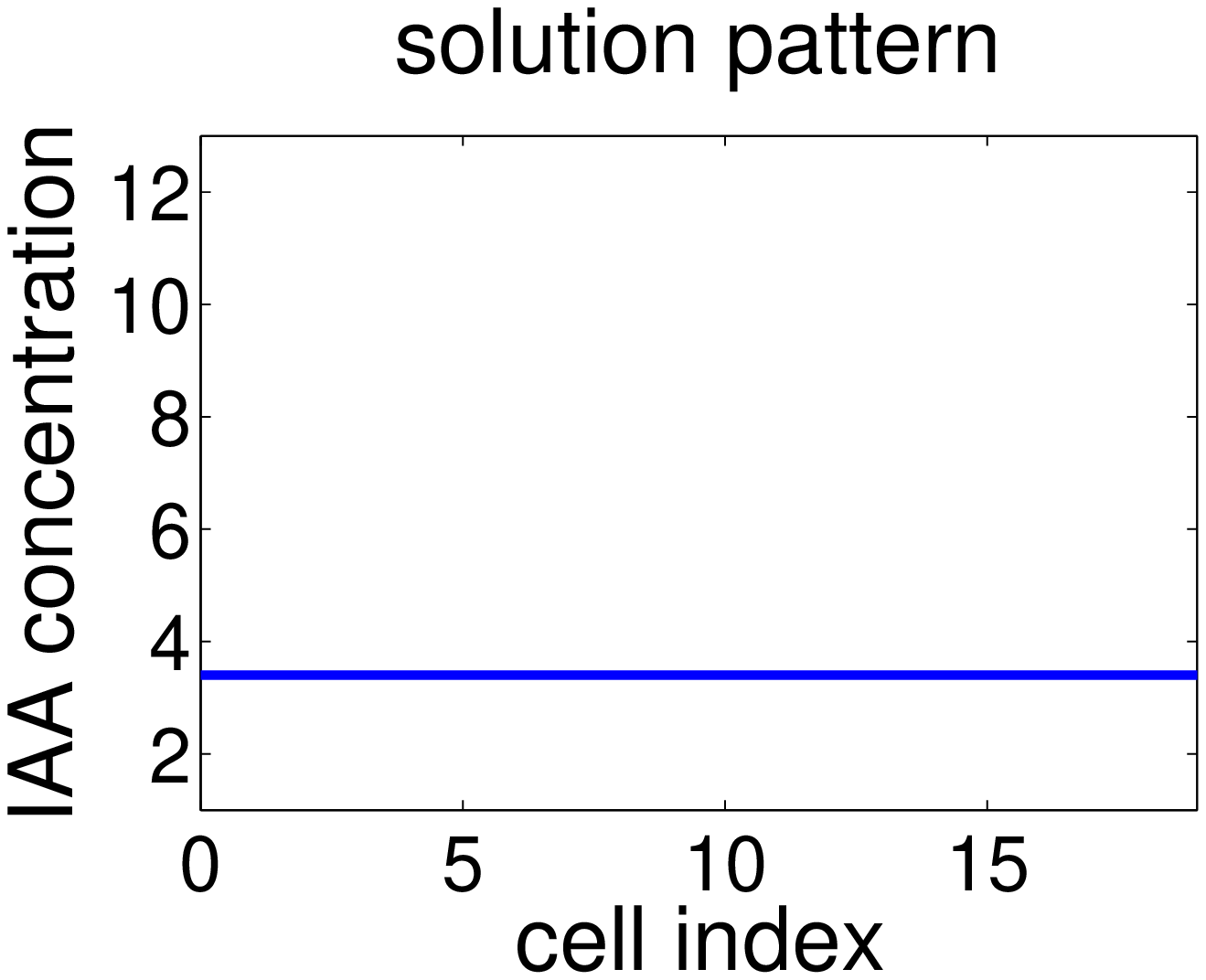}}
	\subfloat[Solution 2] {\label{fig:solution3}\includegraphics[width = 0.3\textwidth]{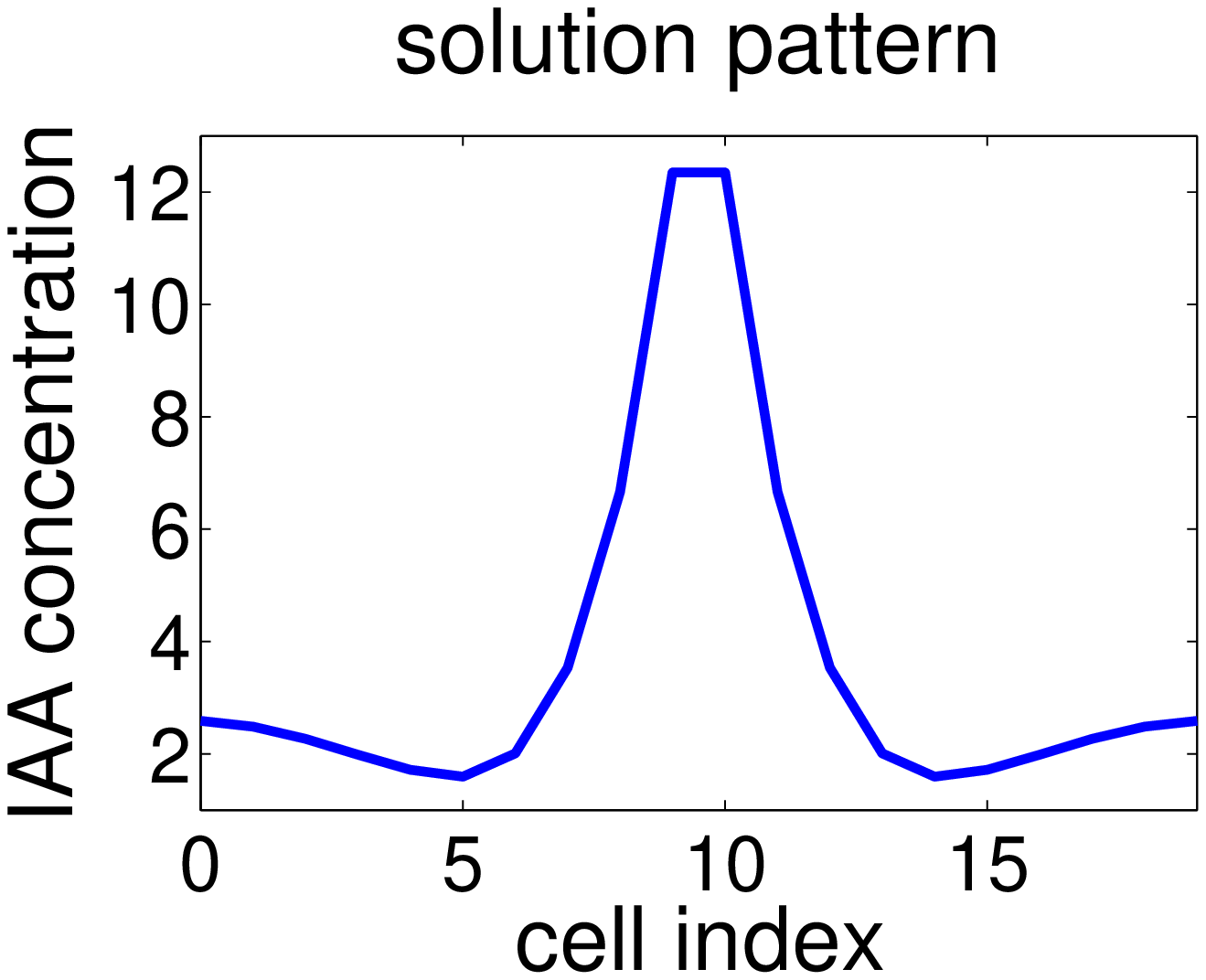}}
	\subfloat[Solution 3] {\label{fig:solution4}\includegraphics[width = 0.3\textwidth]{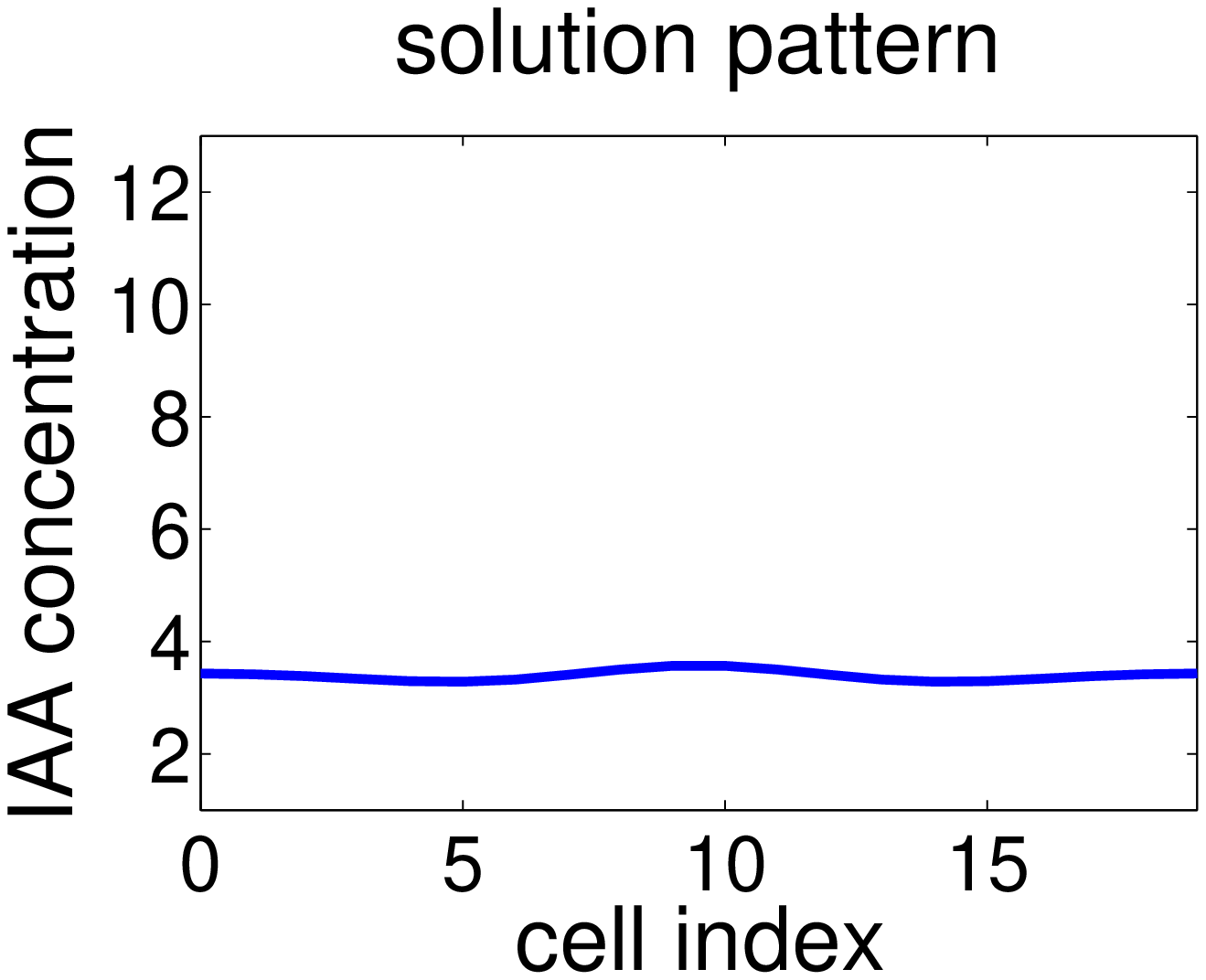}}\\
	\subfloat[Solution 4] {\label{fig:solution5}\includegraphics[width = 0.3\textwidth]{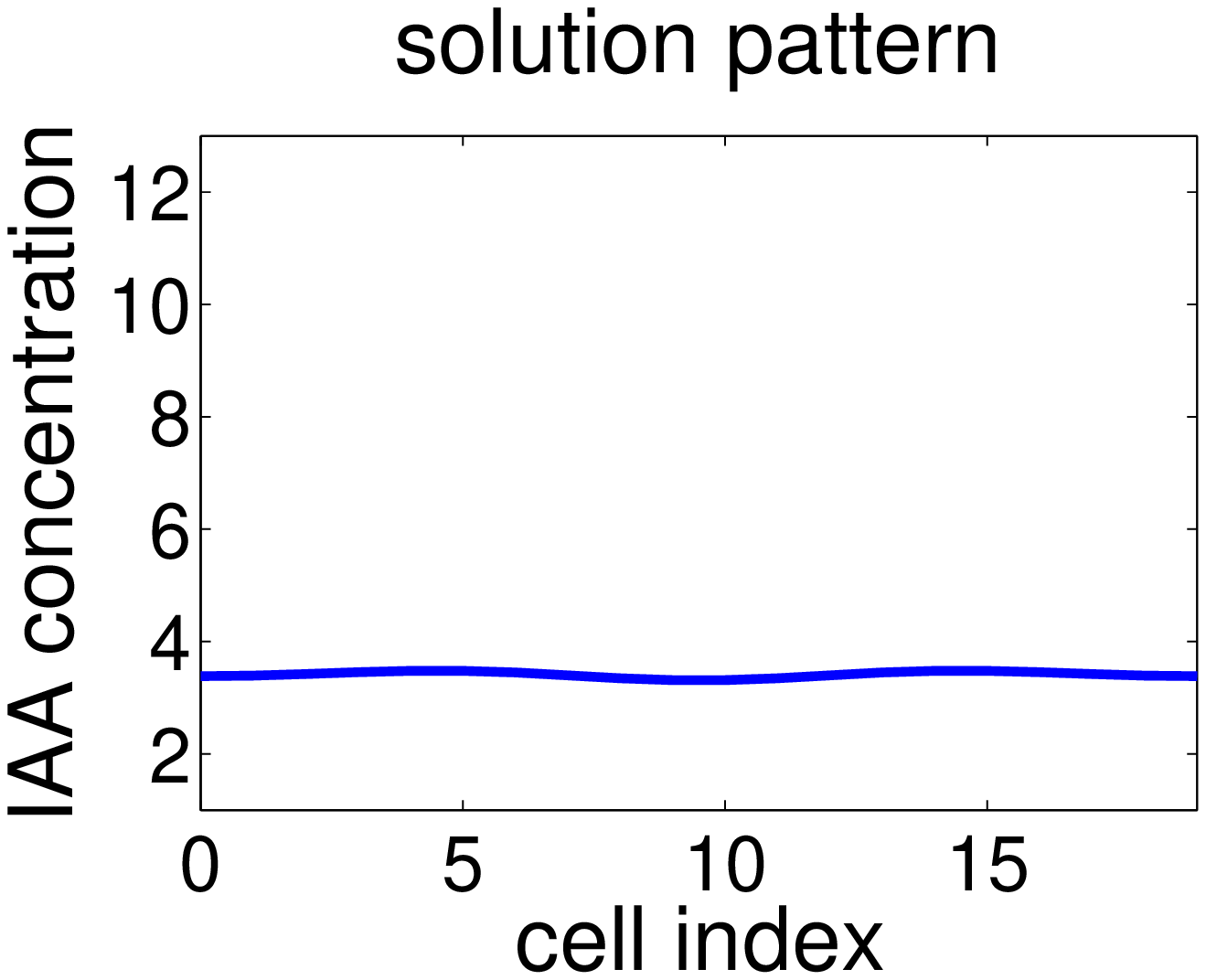}}
	\subfloat[Solution 5] {\label{fig:solution7}\includegraphics[width = 0.3\textwidth]{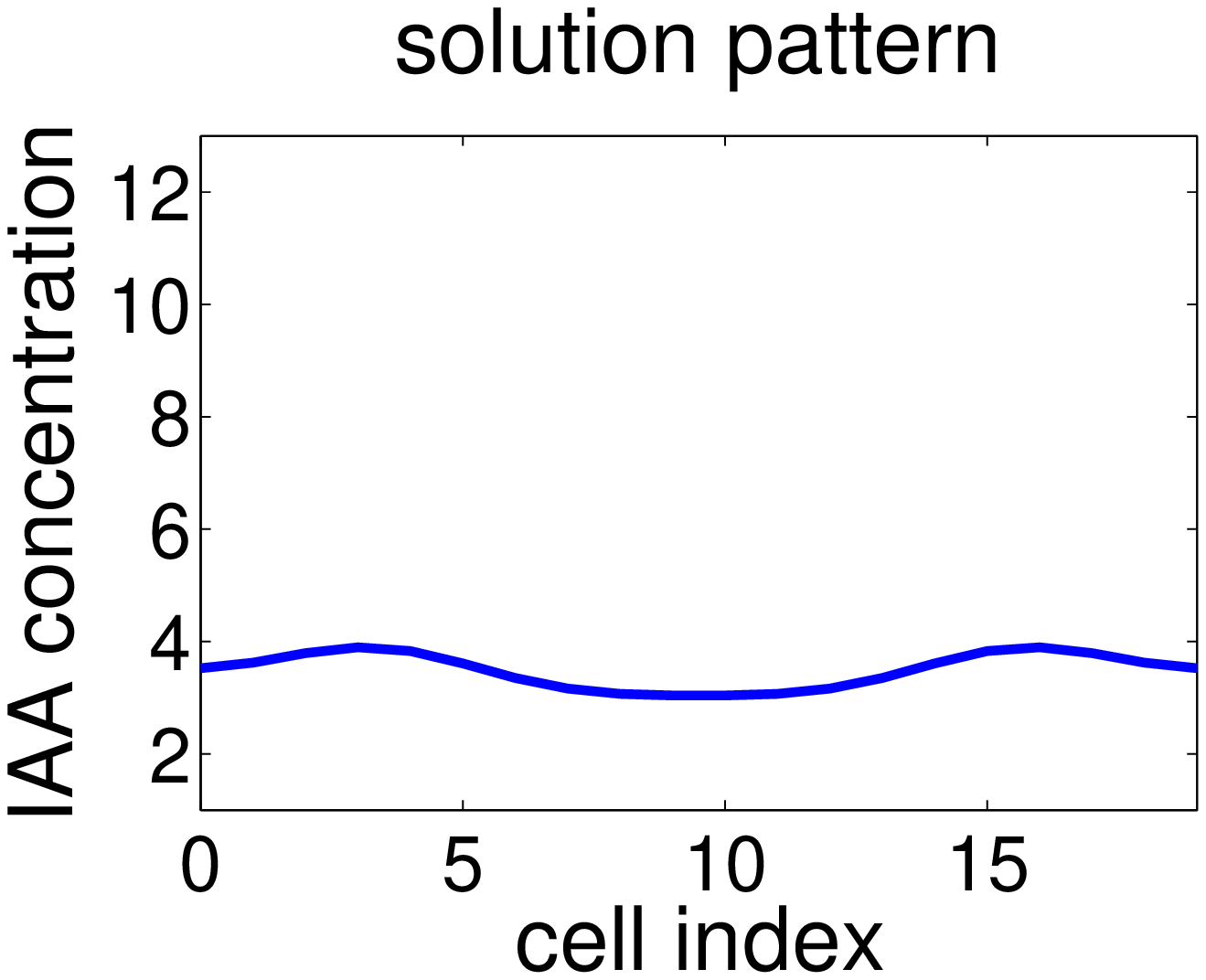}}
	\subfloat[Solution 6] {\label{fig:solution8}\includegraphics[width = 0.3\textwidth]{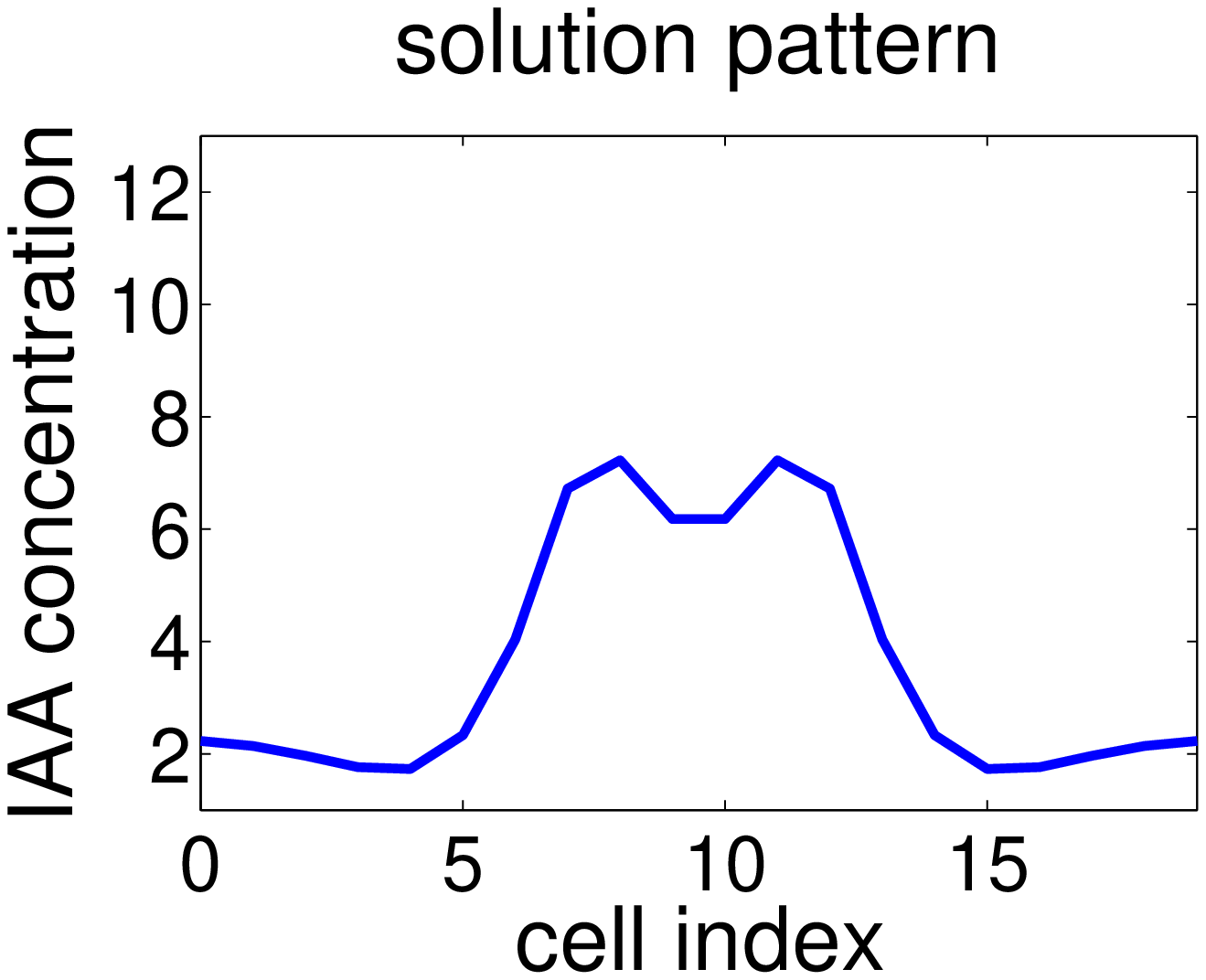}}
	\caption{(a) 
				The bifurcation diagram of example 1 \change{(steady state 
				equation \eqref{eq:Neumannsteadystate})} for a row of 
				$20$ cells with the IAA concentration 
				in cell number 6 versus the continuation parameter $\textrm{T}$ 
				(IAA transport coefficient). Other parameters are taken from M1. 
				BP denotes a branch point and LP a limit point. 
				The stars mark the places of the figures displayed below. \\
        (b)-(g) 
        On these figures the IAA concentration in the whole domain is
        displayed corresponding with the stars marked on
        (a). 
         }
	\label{fig:bifurcationdiagramsolutions}
\end{figure}

\paragraph{Example 1.}  \change{This example illustrates the first
  generic scenario that is encountered when the trivial solution loses
  its stability. The results of the bifurcation analysis for equation \eqref{eq:Neumannsteadystate}
 and parameter set M1 are shown in figure
\ref{fig:bifurcationdiagramsolutions}.}
Figure \ref{fig:BifurcationDiagramScenario1} shows the bifurcation
diagram that depicts the concentration of auxin in
cell number 6 versus the parameter T. The other plots in figure
\ref{fig:bifurcationdiagramsolutions} show the steady state auxin
patterns in all cells for the specific places indicated with
labels in the bifurcation diagram. \\
The trivial solution curve is the starting point of the
continuation. It is the flat horizontal line in the bifurcation
diagram. When the parameter T becomes larger than a critical value
(T$=0.8983$)
, the trivial solution loses its stability at a branch point.  It was
found by calculating for every solution
$\left(\textbf{U}_i,\boldsymbol{\lambda}_i\right)$ on the branch, the eigenvalues
of the augmented Jacobian matrix defined as
\begin{equation}
J_{\text{aug}} = \left[J_{\mathbf{U}} \left| J_{\boldsymbol{\lambda}}\right.\right].
\label{eq:AugmentedJacobian}
\end{equation}
If the Jacobian $J_{\mathbf{U}}$ is singular and the rank of the
augmented Jacobian is still smaller than $2n$, then the solution point
$\left(\textbf{U}_i,\boldsymbol{\lambda}_i\right)$ is a branch point. This means
that there exist an eigenvalue $\mu\left(\boldsymbol{\lambda}_i\right)$ of the
Jacobian which is equal to zero. Inserted into a graph, there is a
path of an eigenvalue of the solution points corresponding to
$\boldsymbol{\lambda}$ close to $\boldsymbol{\lambda}_i$, that crosses the imaginary axis at the
real axis when $\boldsymbol{\lambda} = \boldsymbol{\lambda}_i$.  In the branch point on figure
\ref{fig:bifurcationdiagramsolutions} there is an exchange of
stability to another branch, also shown in the diagram. There are two
stable parts on this other solution branch with patterns.  When the
IAA transport coefficient T is large, the stable solution pattern on
this branch consist of one big peak (figure \ref{fig:solution3}).
The other stable part on this branch appears in a very limited range
where T is smaller. For example solution 4 is such a stable pattern
and it has two small variations (figure \ref{fig:solution5}).  The
pattern in figure \ref{fig:solution3} is the same pattern that was
obtained by numerical integration with the fourth order Runge-Kutta in
figure \ref{fig:timeEvolution} as discussed in Section
\ref{subsec:TimeIntegration}. We thus found a connection between the
trivial flat solution and the numerical solution with peaks.

\paragraph{Example 2.} \change{Here we illustrate a second scenario
  in which the trivial solution loses its stability.}  It describes the
results of the bifurcation analysis for the model with parameter set
M2 that Smith \textit{et al.} used in their publication
\citep{Smith2006}. It differs from the parameter set M1 by a lower
production coefficient of IAA, $\rho_{_{\IAA}}$. In the previous
example, the stability was lost in a branch point.  Now, we find that
the stability is lost through a Hopf bifurcation where the equilibrium
transitions into a periodic orbit.  Looking at the eigenvalues of the
Jacobian in this Hopf point, there is a pair of eigenvalues that
satisfies
\begin{equation}
\mu\left(\boldsymbol{\lambda}_i\right) = \pm i\beta.
\end{equation} 
If we draw a trajectory of the eigenvalues of solution points with
$\boldsymbol{\lambda}$ close to $\boldsymbol{\lambda}_i$, we see that there are two complex
conjugated eigenvalues different from zero that cross the imaginary
axis when $\boldsymbol{\lambda} = \boldsymbol{\lambda}_i$.
\begin{figure}
	\centering
	\subfloat[Bifurcation Diagram] {\label{fig:BifurcationDiagramScenario2} \includegraphics[width = 0.95\textwidth]{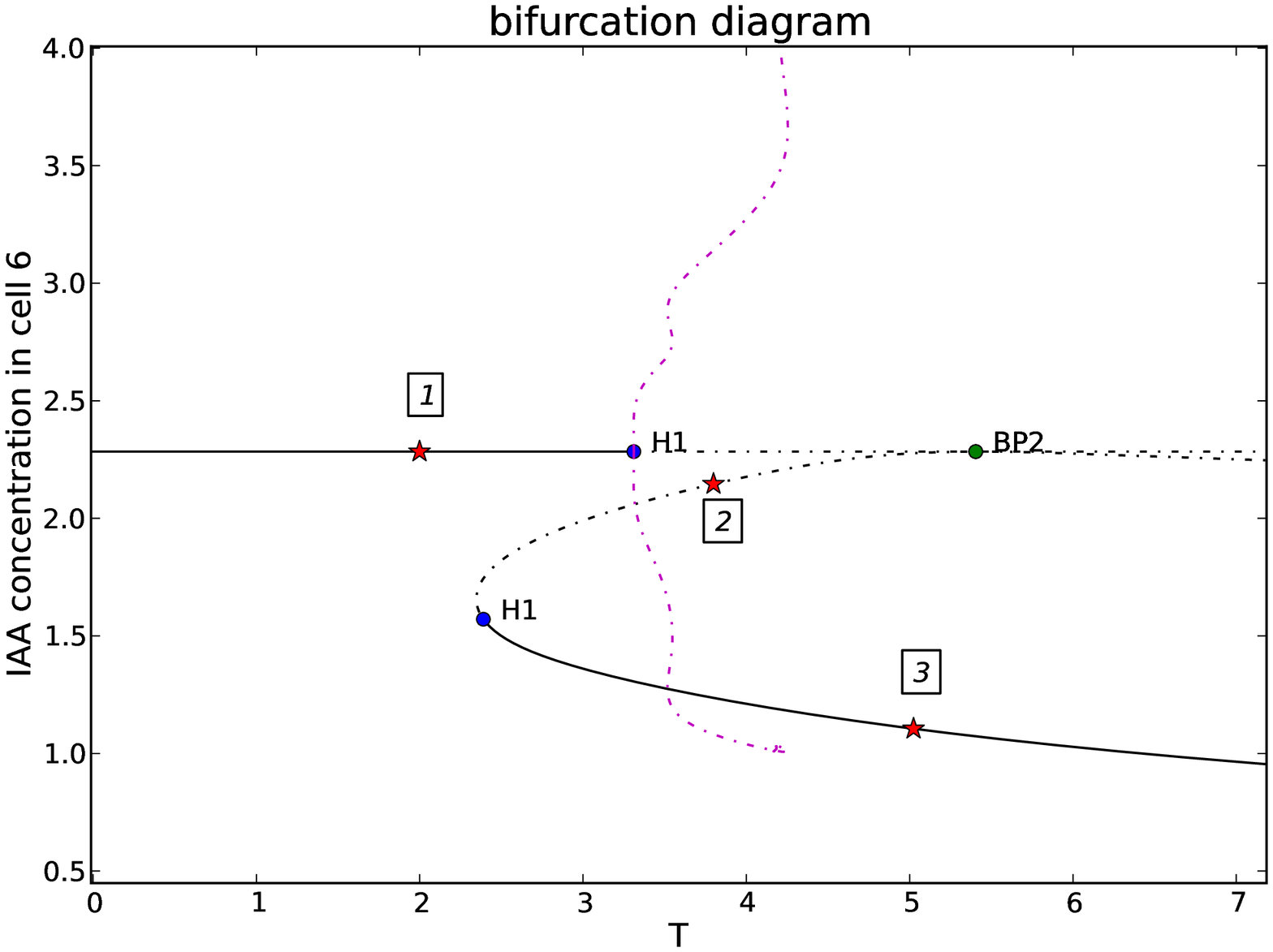}}\\
	\subfloat[Solution 1] {\label{fig:solution1hopf}\includegraphics[width = 0.3\textwidth]{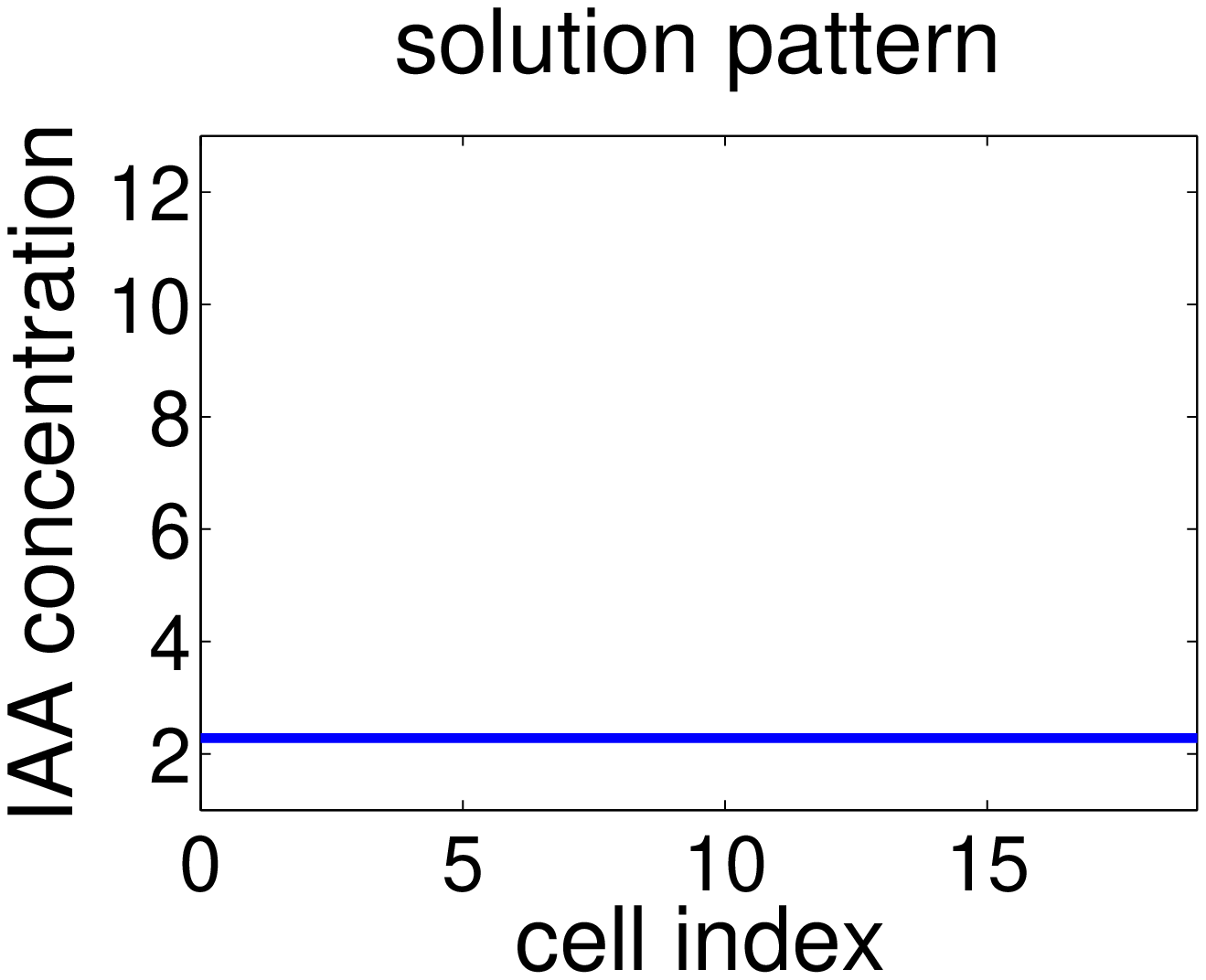}}
	\subfloat[Solution 2] {\label{fig:solution3hopf}\includegraphics[width = 0.3\textwidth]{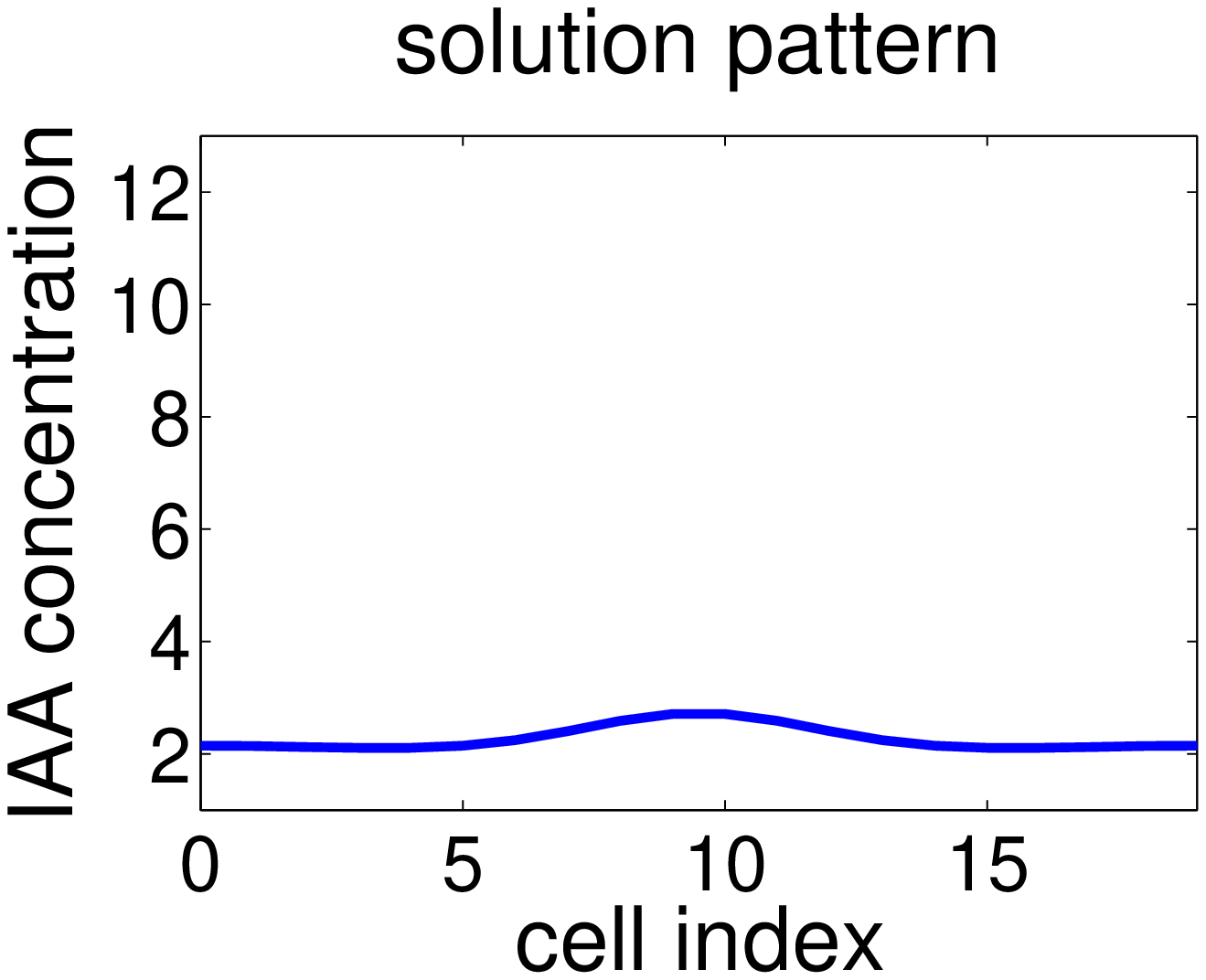}}
	\subfloat[Solution 3] {\label{fig:solution4hopf}\includegraphics[width = 0.3\textwidth]{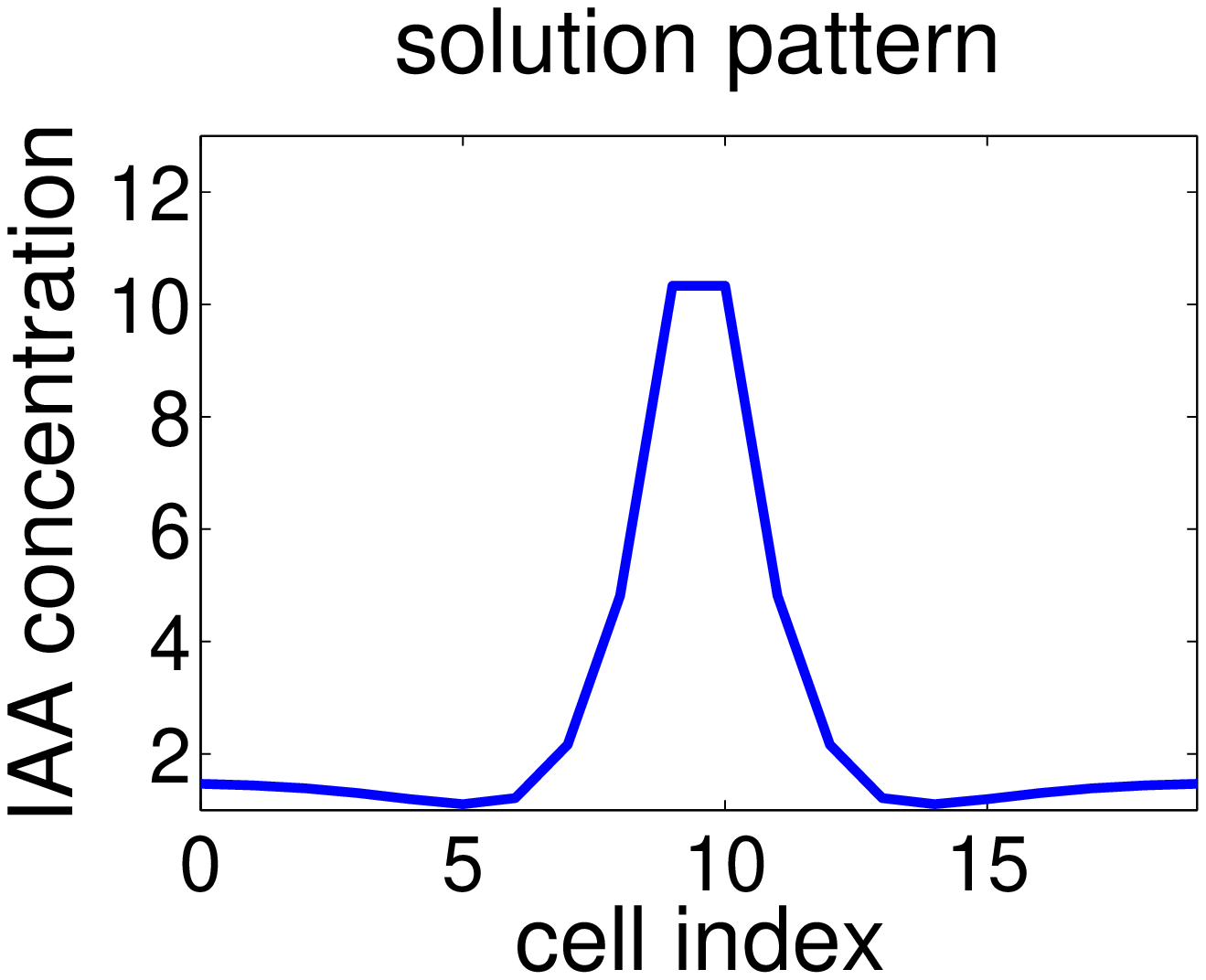}}
	\caption{(a) 
						The bifurcation diagram of example 2 \change{(steady state equation \eqref{eq:Neumannsteadystate})}
						for a row of $20$ cells
          	with the IAA concentration in cell number $6$ versus the 
          	continuation parameter $\textrm{T}$ (IAA transport
          	coefficient). Other parameters are taken from M2. The stars mark
          	the places of the figures displayed below. The dotted line 
          	through point H1 shows the maximal and minimal value of the 
          	IAA concentration in cell number 6 of the periodic solution 
          	for each choice of the parameter T. H and BP denote respectively 
          	a Hopf point and a branch point.\\
            (b)-(d) 
            On these figures the IAA concentration in the whole domain is
          	displayed corresponding with the stars marked on
          	(a). 
          	  }
	\label{fig:bifurcationdiagramsolutionshopf}
\end{figure}

Figure \ref{fig:BifurcationDiagramScenario2} shows the bifurcation
diagram depicting again the concentration of IAA in cell $6$ versus
the continuation parameter T \change{for the same model as in example 1,} 
but now for parameter set M2. In this
situation, the stability of the trivial solution is lost in a Hopf
point at T $=3.3113$. 
The branch that emerges from this Hopf point, shows the maximal and
minimal IAA concentration over the orbit for each choice of parameter
T.  All the solutions on this branch are unstable and therefore we only
have unstable periodic solutions. Further, also another steady state
branch, different from the trivial solution branch, is displayed.
This branch intersects with the trivial solution branch at a branch
point (T $=5.4047$ 
). Around this branch point, all solutions are unstable. However, when
we follow this new branch, we encounter another Hopf point where we
now gain stability. The pattern of these stable solutions consist of
one single big peak in the middle of the domain (see figure
\ref{fig:solution4hopf}).

\paragraph{Example 3.} \change{This example illustrates that there are
  stable orbits beyond the Hopf bifurcation point for some particular
  choices of the parameters.} For the third example parameter set M3 is
used that differs from sets M1 and M2 in the production coefficient
of IAA. It is smaller than in set M2.
\change{Again equation \eqref{eq:Neumannsteadystate} is solved.}

\begin{figure}
	\centering
	\includegraphics[width = 0.95\textwidth]{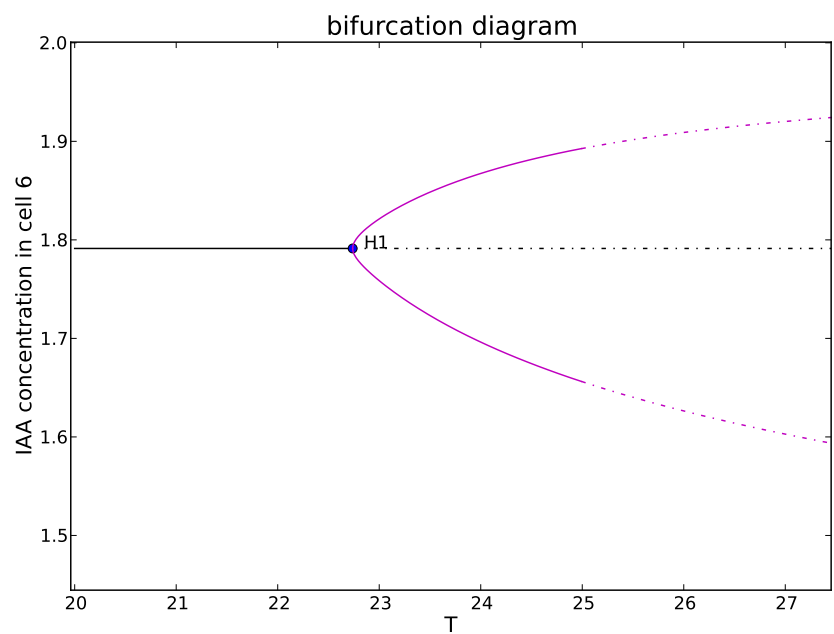}
	\caption{Bifurcation diagram of example 3 \change{where equation \eqref{eq:Neumannsteadystate} was solved} with parameter set M3. It depicts the IAA concentration in cell number 6 versus the continuation parameter $\textrm{T}$ (IAA transport coefficient). H denotes a Hopf point. The stable orbit for T$=23.5$ is shown on figure \ref{fig:TimeEvolutionPeriodic}.}
	\label{fig:bifurcationdiagramM3} 
\end{figure}

\begin{figure}
	\centering
	\includegraphics[width = 0.95\textwidth]{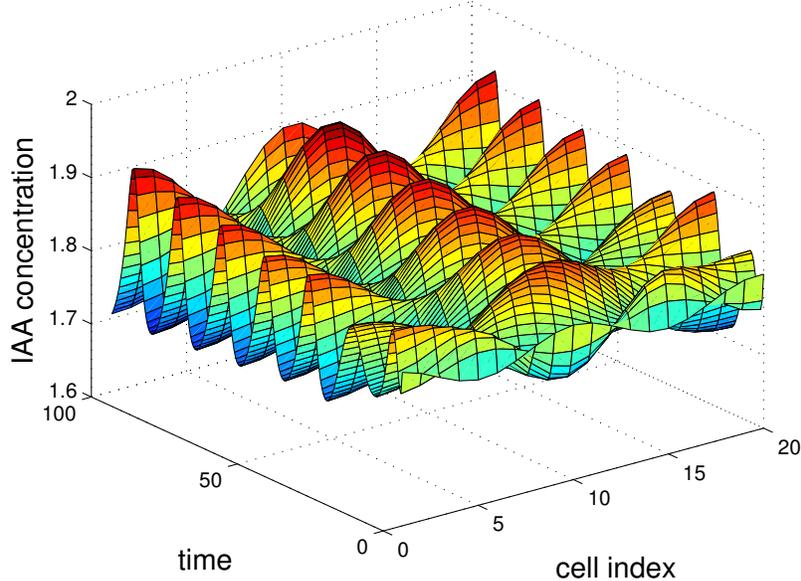}
	\caption{The time evolution \change{of equations \eqref{eq:PIN_p}, \eqref{eq:IAA_a} and \eqref{eq:ActiveTransport}} for a row of $20$ cells starting from the initial value in equation \eqref{eq:beginvalue2} and with parameter set M3 but with IAA transport coefficient T$=23.5$. We used RK4 for numerical integration. The resulting solution is a periodic solution where the pattern changes from one peak concentration of auxin in the middle of the domain to a pattern with high concentrations at the boundaries. The periodic solution corresponds with the solution for T$=23.5$ on the periodic solution branch in figure \ref{fig:bifurcationdiagramM3}.}
	\label{fig:TimeEvolutionPeriodic}
\end{figure}

The resulting bifurcation diagram is shown in figure
\ref{fig:bifurcationdiagramM3}.
As in example 2, the stability of the trivial solution is lost
in a Hopf point (T $=22.7384$).  
However, in contrast with this previous example, the periodic solution
branch that intersects with the trivial solution branch in this point
contains stable periodic solutions. Figure
\ref{fig:TimeEvolutionPeriodic} shows, in a three dimensional plot,
the stable periodic solution for IAA transport coefficient $T=23.5$
found with RK4 starting from the initial value
\begin{equation}
\label{eq:beginvalue2}
p_i(0) = 1.79 \quad \text{and} \quad   a_i(0)=3.76,
\end{equation}
where a small perturbation $0.02\sin\left(\left(5\left(i+2\right)\pi\right)/24\right)$ for $i=1,\ldots,20$ was added.
\change{Any other initial state nearby will lead to the same long term solution.}
We see that the periodic solution changes in time from a pattern with
one peak concentration of auxin in the middle of the domain to a
pattern with two high auxin concentrations at the sides of the
domain.
\begin{figure}
	\centering
	\includegraphics[width = 0.95\textwidth]{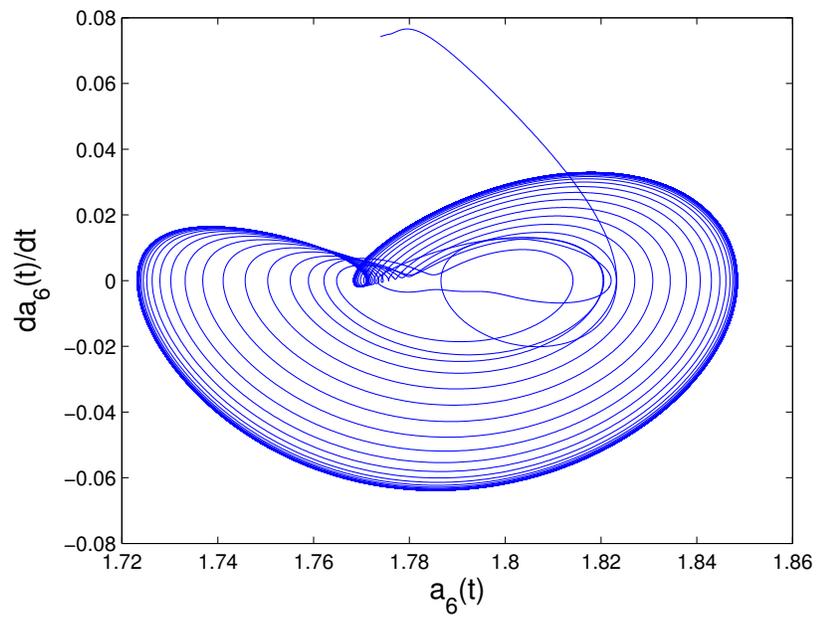}
	\caption{The trajectory of the time evolution of figure \ref{fig:TimeEvolutionPeriodic} in the $\left(a_6\left(t\right), da_6\left(t\right)/dt\right)$-plane. We used a row of $20$ cells starting from initial value \eqref{eq:beginvalue2} and parameter set M3 but with IAA transport coefficient T$=23.5$ and we used RK4 for numerical integration. }
	\label{fig:traject} 
\end{figure}
In figure \ref{fig:traject} we plotted this trajectory in the $\left(a_6\left(t\right), da_6\left(t\right)/dt\right)$-plane starting from the initial value in equation \eqref{eq:beginvalue2}.

\begin{figure}
	\centering
	\subfloat[Bifurcation Diagram] {\label{fig:BifurcationDiagram100} \includegraphics[width = 0.95\textwidth]{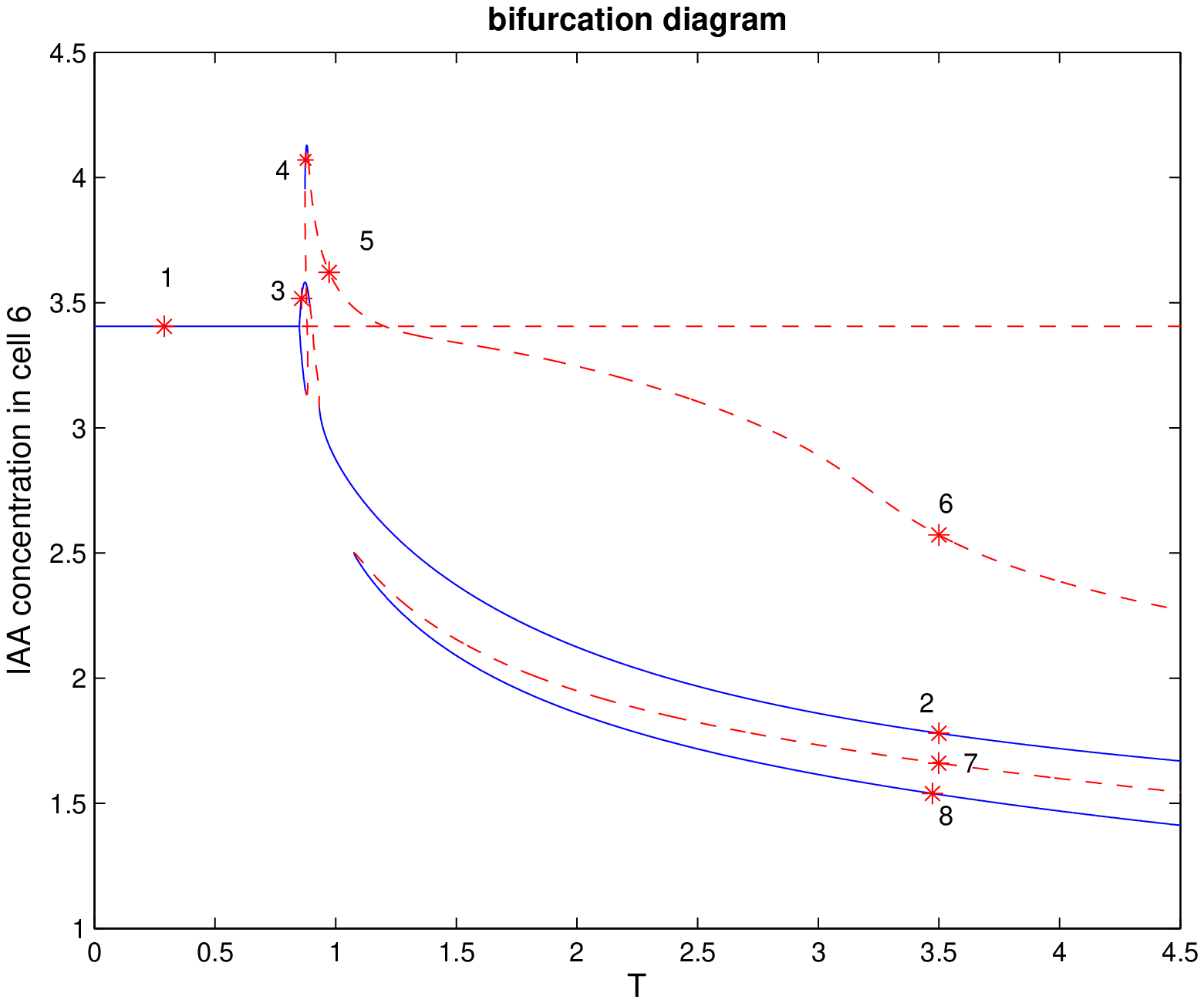}}\\
	\subfloat[Solution 1] {\label{fig:solution1_100}\includegraphics[width = 0.24\textwidth]{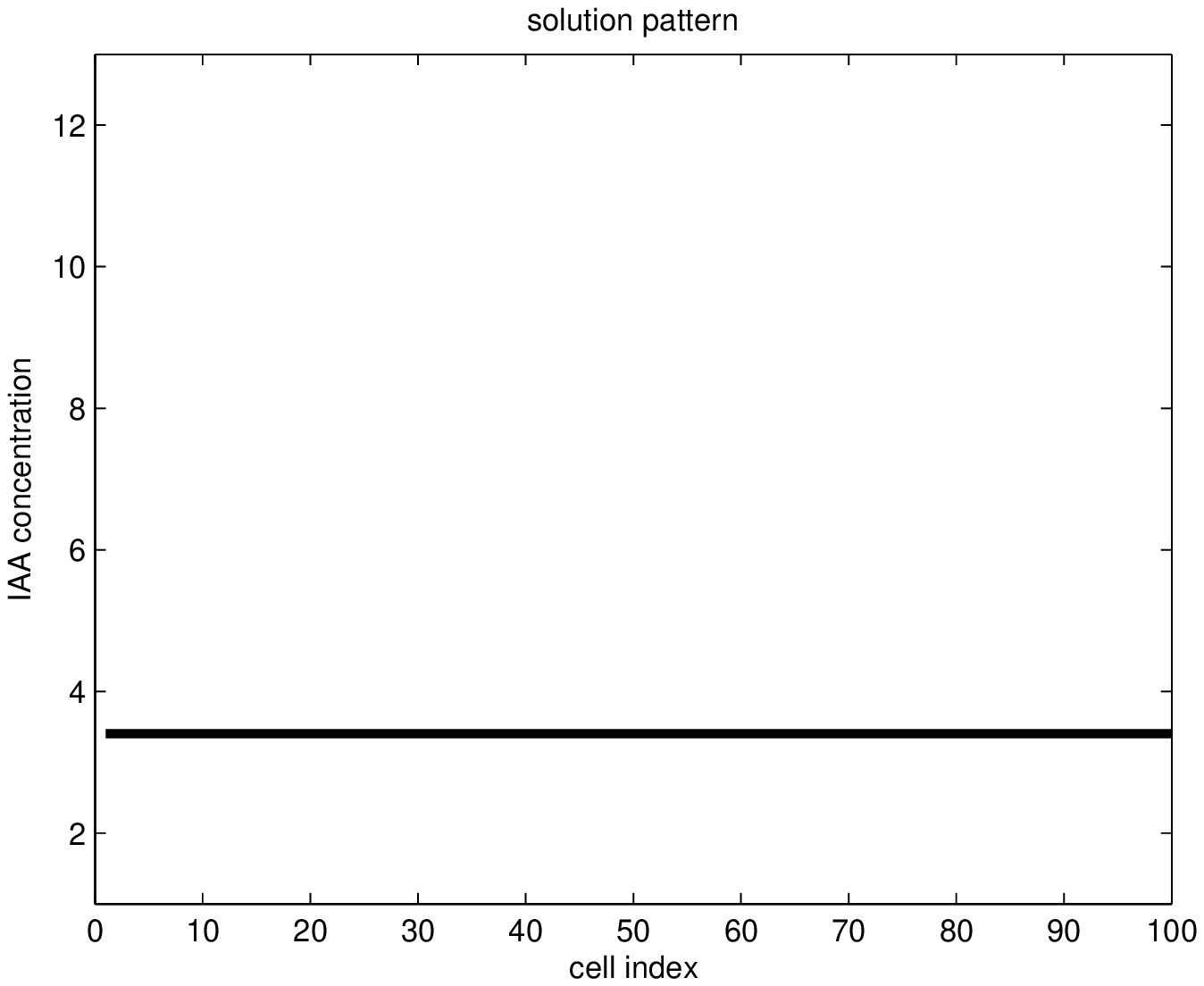}}
	\subfloat[Solution 2] {\label{fig:solution2_100}\includegraphics[width = 0.24\textwidth]{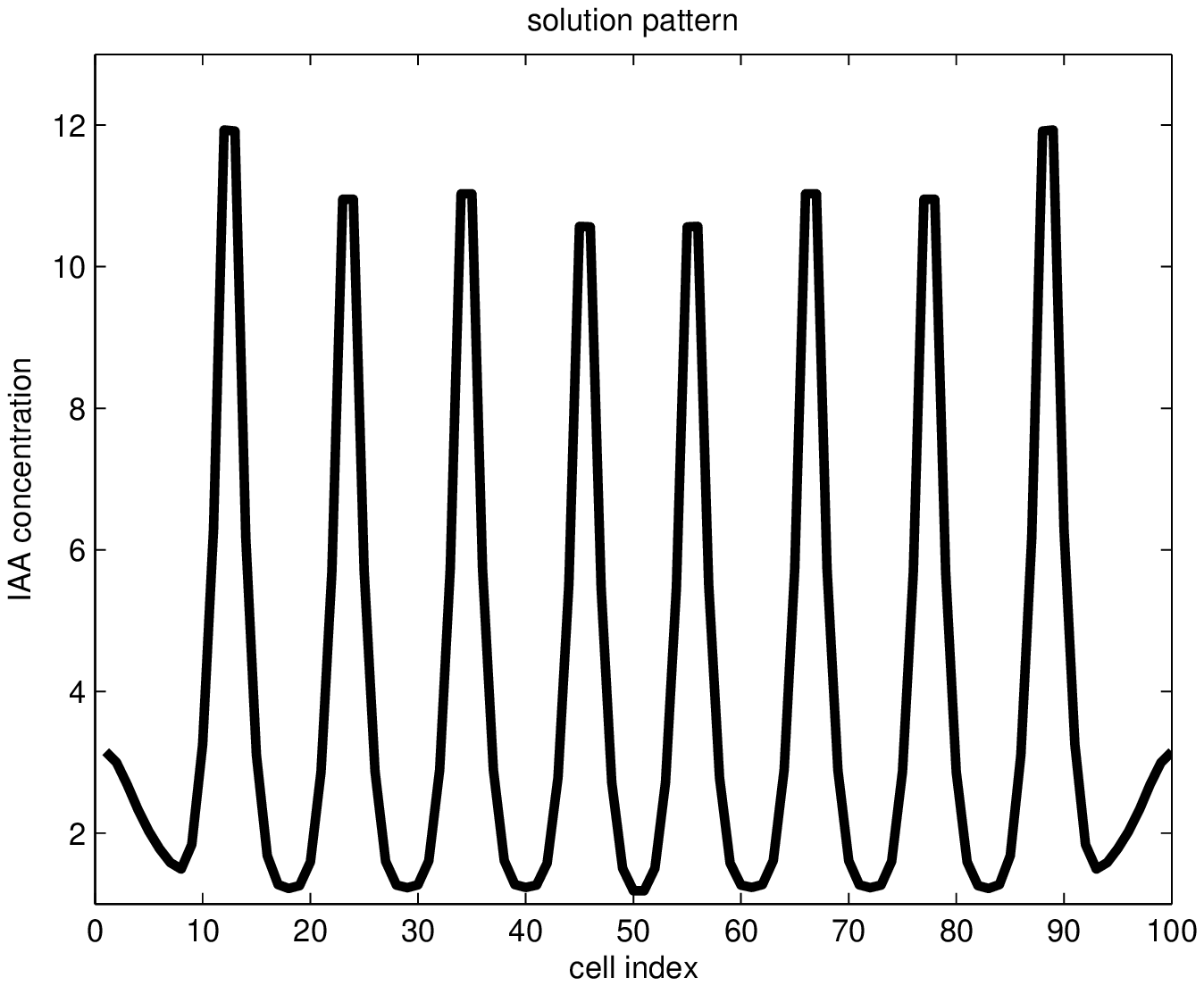}}
	\subfloat[Solution 3] {\label{fig:solution3_100}\includegraphics[width = 0.24\textwidth]{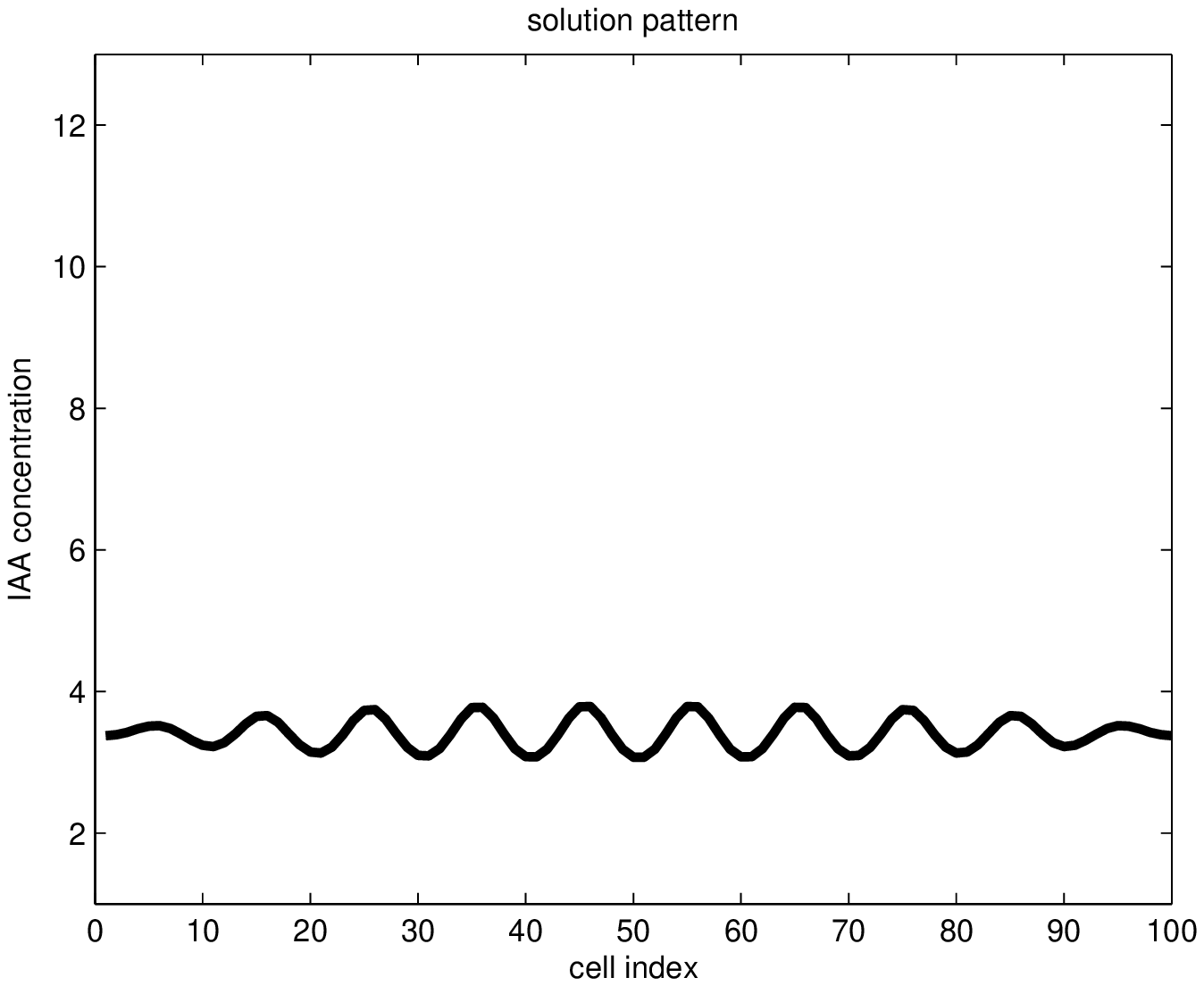}}
	\subfloat[Solution 4] {\label{fig:solution4_100}\includegraphics[width = 0.24\textwidth]{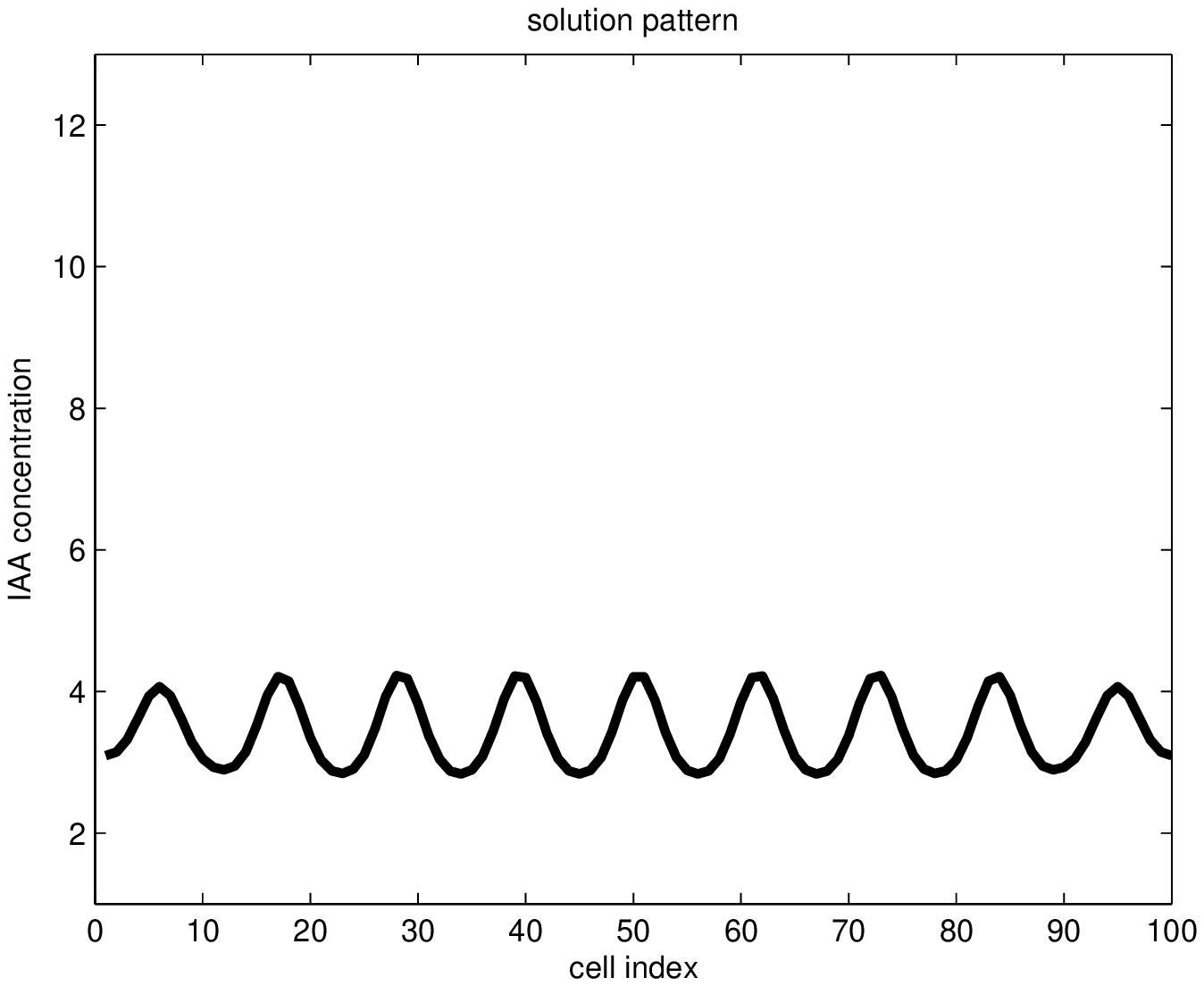}} \\
	\subfloat[Solution 5] {\label{fig:solution5_100}\includegraphics[width = 0.24\textwidth]{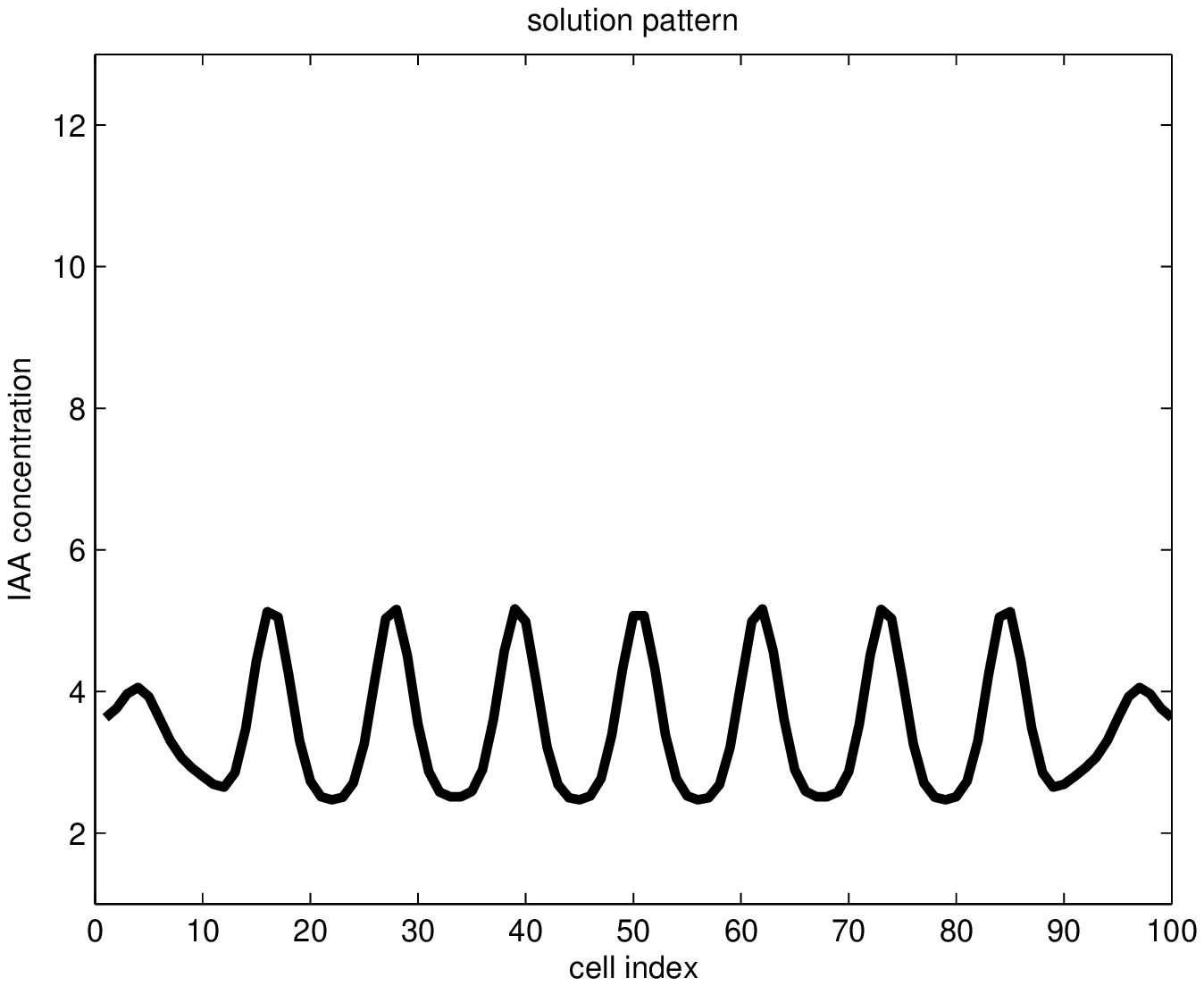}}
	\subfloat[Solution 6] {\label{fig:solution6_100}\includegraphics[width = 0.24\textwidth]{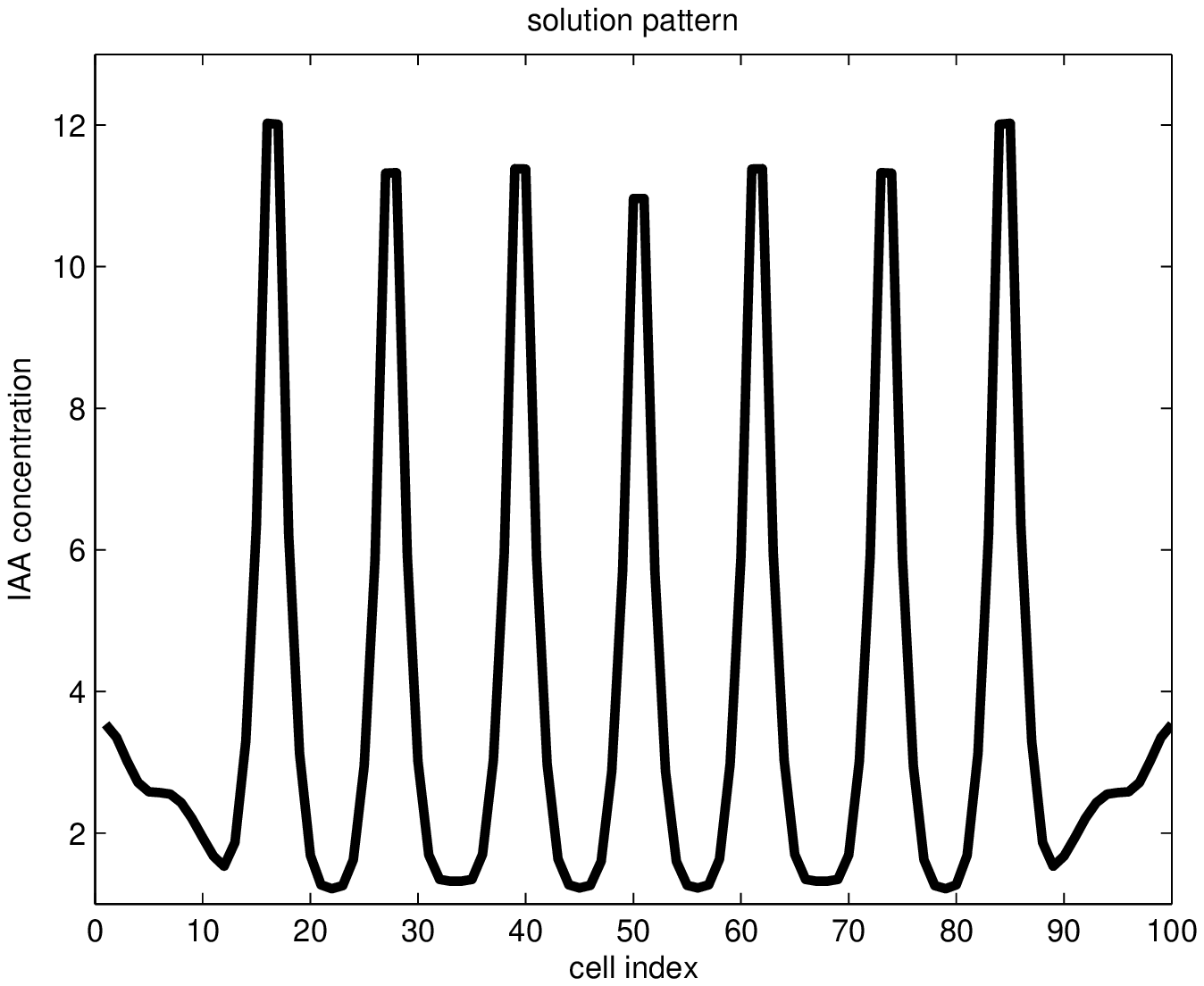}}
	\subfloat[Solution 7] {\label{fig:solution7_100}\includegraphics[width = 0.24\textwidth]{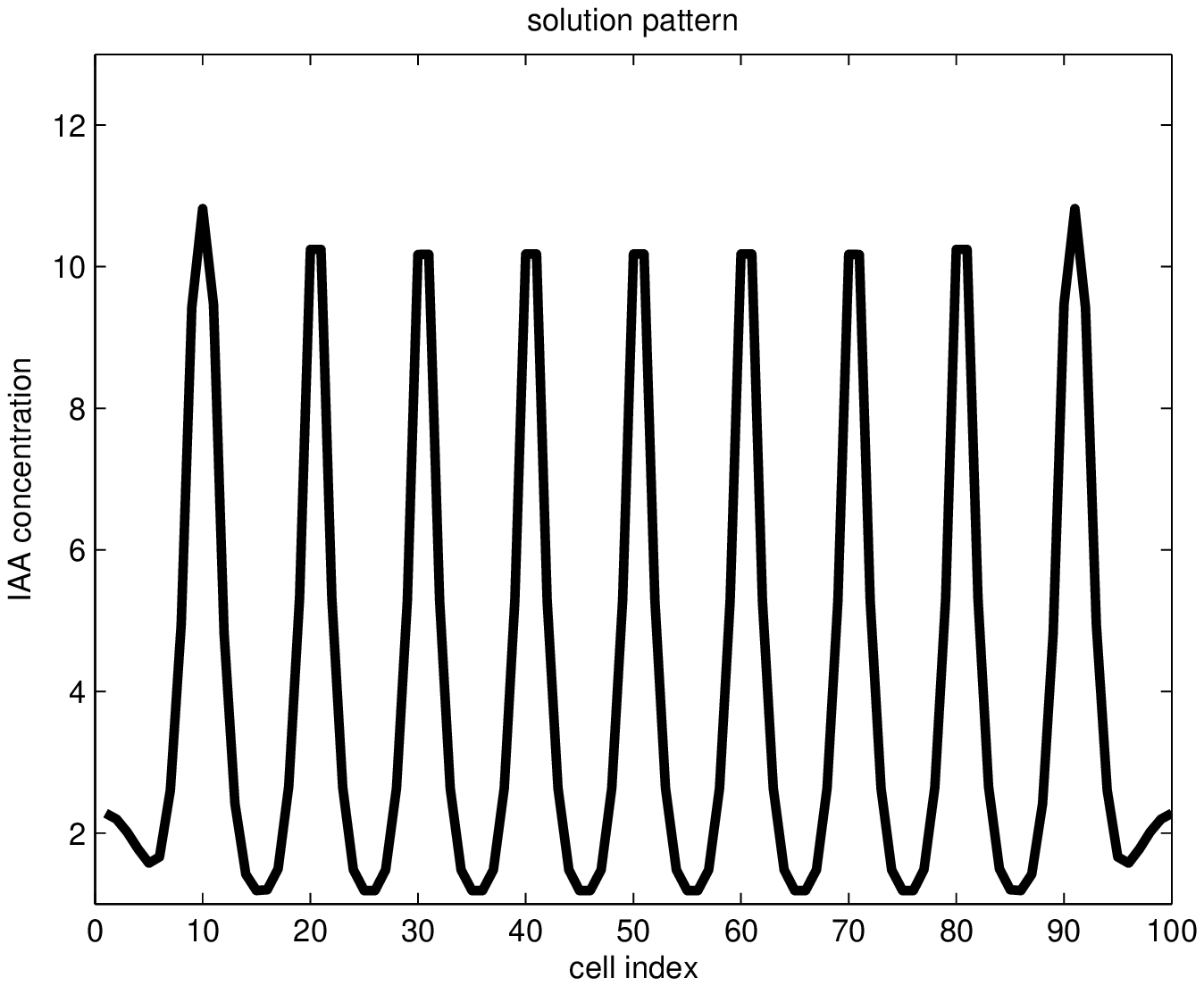}}
	\subfloat[Solution 8] {\label{fig:solution8_100}\includegraphics[width = 0.24\textwidth]{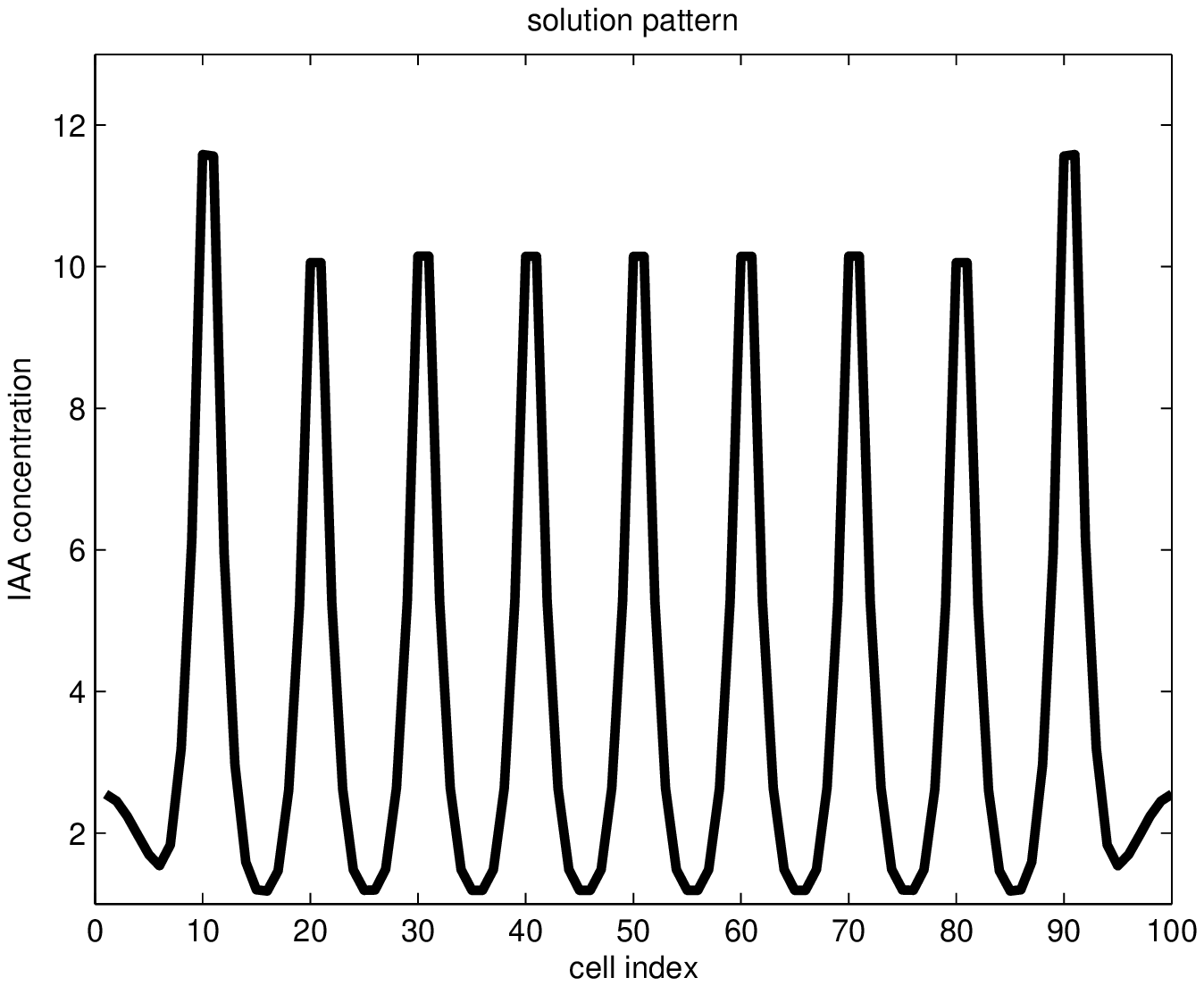}}
	\caption{(a) 
					Bifurcation
          diagram of example 4 \change{(steady state equation \eqref{eq:Neumannsteadystate})}
          for a row of $100$ cells and parameter
          set M1.  The diagram shows the IAA
          concentration in cell number $6$ versus the continuation parameter
          $\textrm{T}$ (IAA transport coefficient). The unconnected branch is
          found by time integrating the steady state solution of
          figure \ref{fig:solution3} copied five times in a row until
          a steady state is found. The resulting steady state is then
          used as a starting point for the continuation.The stars
          mark the places of the figures displayed below.
          (b)-(i) 
          On
          these figures the IAA concentration in the whole domain is
          displayed corresponding with the stars marked on
          (a). 
            }
	\label{fig:bifurcationdiagramsolutions100}
\end{figure}

\paragraph{Example 4.} The previous three examples showed a part of
the bifurcation diagrams \change{for equation
  \eqref{eq:Neumannsteadystate}} corresponding with parameter sets M1,
M2 and M3, that differ in IAA production rate, for a one dimensional
domain of 20 cells. In this example we look to a row of $100$ cells
and use the parameter values of set M1. In this example, the stability
of the trivial solution is again, as in example 1, lost at a branch
point (T $=0.8504$) (see figure \ref{fig:BifurcationDiagram100}). The
branch that crosses the trivial solution branch in this point is a bit
more complicated. The branch contains 3 different stable parts.
The stable part that contains solution point 2 consists of solutions
with a pattern with 8 peaks (see figure \ref{fig:solution2_100}). Also
the small stable area on the branch that contains solution 3 consists
of patterns with $8$ peaks, but they are smaller due to the small
value of T in this region. We see that the peaks become higher for an
increasing IAA transport coefficient T (compare for example the
patterns in figures
\ref{fig:solution2_100} and
\ref{fig:solution3_100} or in figures \ref{fig:solution5_100}
and \ref{fig:solution6_100}). The third stable part on this
solution branch, appears also in a limited range of T. The patterns
show only $7$ peaks of high auxin concentration (see figure
\ref{fig:solution4_100}). The patterns on the unstable part of the
branch that contains solutions $5$ and $6$ also consist of 7 peaks
(see figure \ref{fig:solution5_100} and \ref{fig:solution6_100})

The other solution branch in figure \ref{fig:BifurcationDiagram100} is
found by time integration by starting with the steady state solution
of figure \ref{fig:timeEvolution} or the pattern in figure
\ref{fig:solution3} copied five times in a row which rapidly leads to
a steady state that can be used as a starting point for the
continuation. We see that this branch is not (directly) linked with
the trivial solution branch and consists of one stable and one
unstable part.  On both parts the patterns consist of $9$ high peaks
of auxin concentration (see figures \ref{fig:solution7_100} and
\ref{fig:solution8_100}).

\begin{figure}
	\centering
	\subfloat[$\textrm{D}=1.000$] {\label{fig:TypeBifPoint_rhoIAA_T} 																																																											\includegraphics[width=0.5\textwidth]{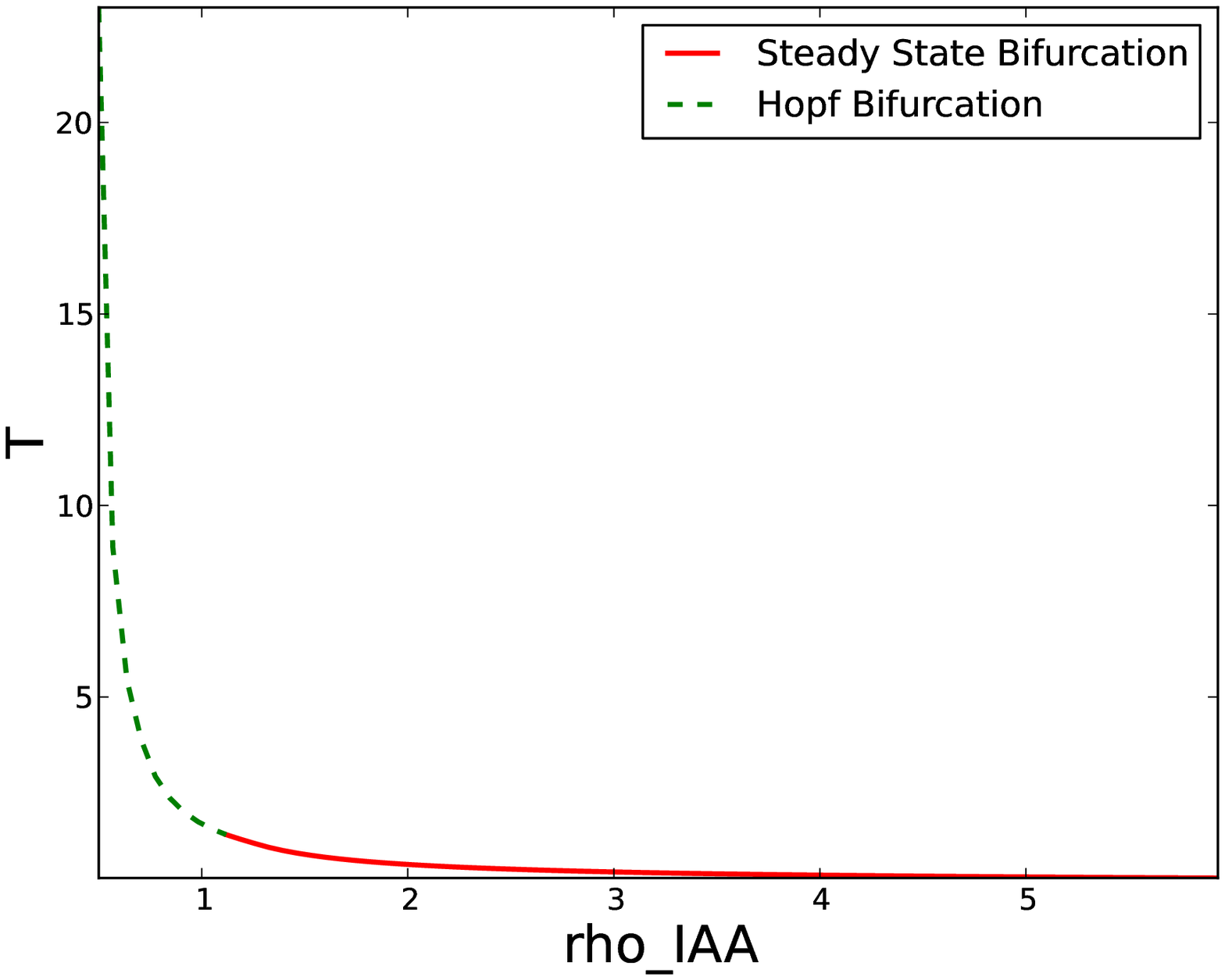}}
	\subfloat[$\textrm{T}=3.500$] {\label{fig:TypeBifPoint_rhoIAA_D} 																																																											\includegraphics[width=0.5\textwidth]{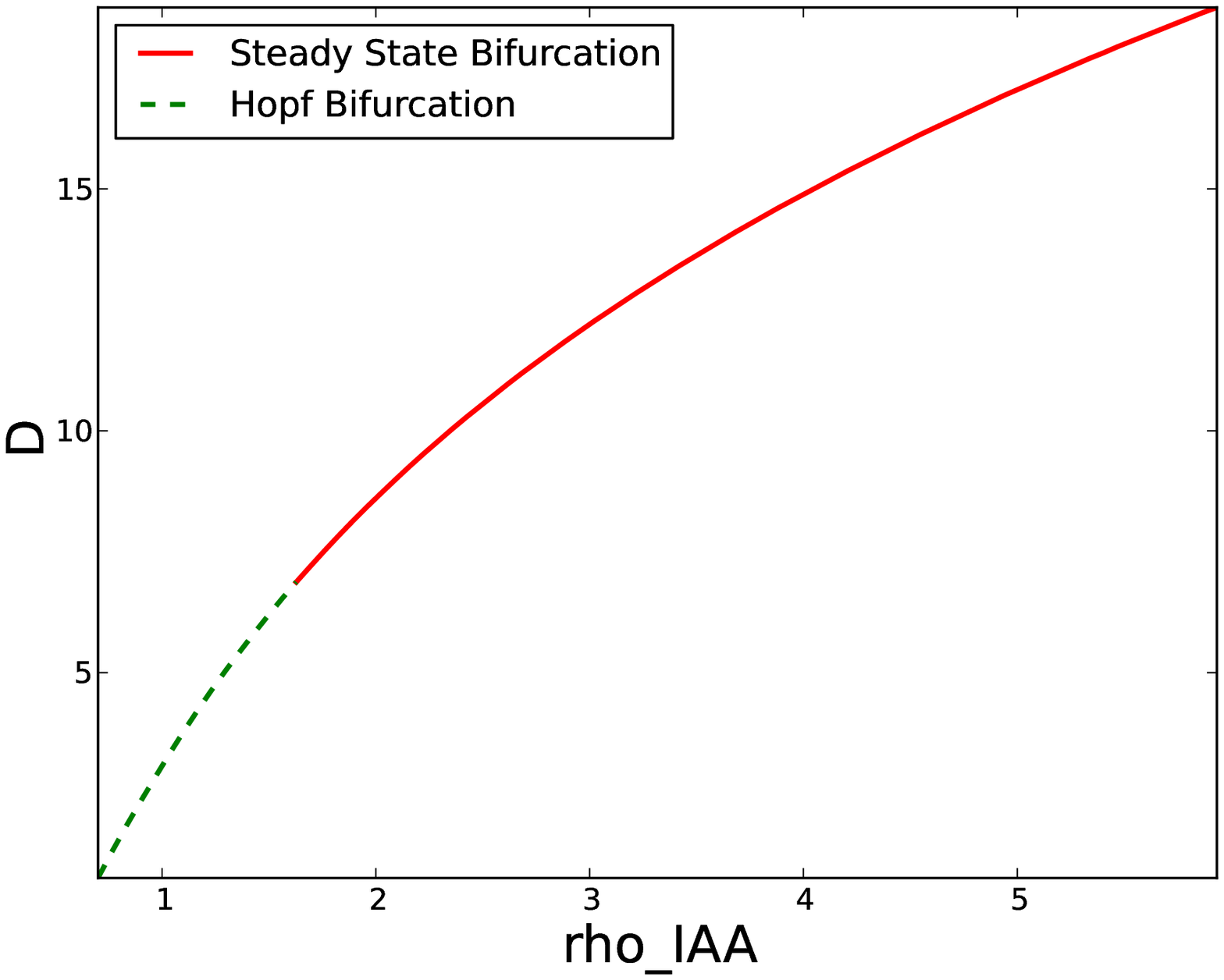}} \\
	\subfloat[$\rho_{_{\IAA}}=0.750$] {\label{fig:TypeBifPoint_D_T}																																																												\includegraphics[width=0.5\textwidth]{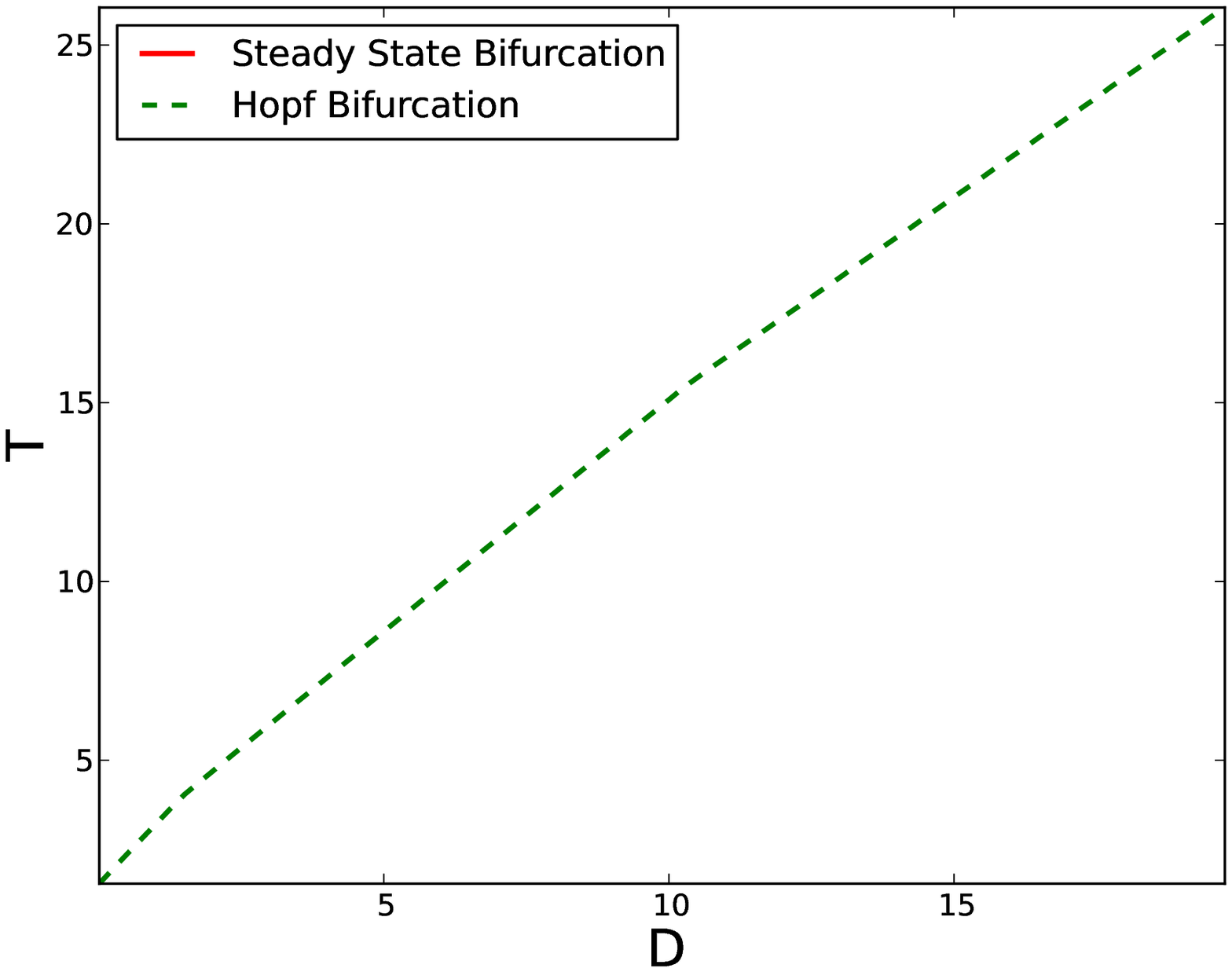}}
	\caption{Map of types of bifurcation points \change{for equation 
					\eqref{eq:Neumannsteadystate} and } a row of $20$
          cells and different choices of the parameters
          $\rho_{_{\IAA}}$ (IAA production coefficient), D (IAA
          diffusion coefficient), T (IAA transport coefficient). Other
          parameters are taken from table
          \ref{table:parametervalues}.}
	\label{fig:TypeBifPoint}
\end{figure}

\paragraph{} Figures \ref{fig:bifurcationdiagramsolutions},
\ref{fig:bifurcationdiagramsolutionshopf},
\ref{fig:bifurcationdiagramM3} and
\ref{fig:bifurcationdiagramsolutions100} for examples 1, 2, 3 and 4
show different bifurcation diagrams \change{of equation
  \eqref{eq:Neumannsteadystate}}. In the examples 1 and 4, where
parameter set M1 is used, the stability of the trivial solution is
lost in a branch point.  While in the examples 2 and 3 it is lost in a
Hopf point. We can calculate the type of these bifurcation points for
every transition from stable to unstable for the three dimensional
space of the parameters D, T and $\rho_{_{\IAA}}$.  We found that in
the one-dimensional case, for a row of cells, only two different
situations can occur: either the stability of the trivial solution is
lost in a Hopf point or it is lost in a branch point. Figure
\ref{fig:TypeBifPoint} shows, for a row of $20$ cells and parameter
values from table \ref{table:parametervalues}, the stability boundary
for the parameters D, T and $\rho_{_{\IAA}}$.  The line also indicates
the corresponding type of bifurcation.
For example figure \ref{fig:TypeBifPoint_rhoIAA_T} shows that for small values of the production coefficient of IAA the stability of the trivial solution will be lost in a Hopf point. Indeed, in example 2 and 3 we found similar results.

\begin{figure}
	\centering
	\includegraphics[height=0.9\textwidth,angle=-90]{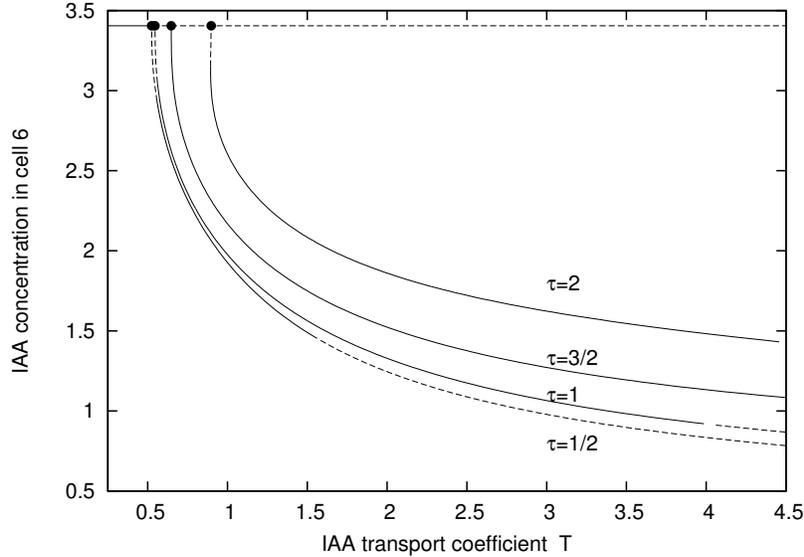}
	\caption{ \change{The bifurcation diagrams of example 5 for a row of
          $20$ cells with the IAA concentration in cell number
          $6$ versus the continuation parameter T (IAA transport
          coefficient). The adapted model with active transport equation
          \eqref{eq:ActiveTransport_all} is used where the exponent in
          the active transport ($\tau$) is set to different values
          ($\tau = 2, 3/2, 1$ and $1/2$). The parameter $\omega$ is
          set to 1 and the other parameter values are taken from M1.}}
	\label{fig:BifDiagram_tau}
\end{figure}

\paragraph{Example 5.} \change{This example explores the role of the
  quadratic dependence of the active transport on the auxin
  concentration and the bifurcation diagrams will be calculated for
  varying exponents $\tau$. We will use the adapted model with active
  transport equation \eqref{eq:ActiveTransport_all} with parameter
  $\omega$ equal to one and for different values for parameter
  $\tau$.  The results in figure \ref{fig:BifDiagram_tau} show the
  concentration of auxin in cell number $6$ versus the parameter T for
  $\tau = 1/2,1,3/2$ and $2$. In this example we use again parameter
  set M1 and a row of $20$ cells, as in example $1$.  }
\change{ The flat branch, on top of the figure, is the trivial
  solution branch, which is independent of the value for the exponent
  $\tau$. However, the stability of this branch depends on $\tau$
  and is only shown for $\tau = 1/2$.  }

\change{ Similarly to example 1 the trivial solution loses, for each
  choice of $\tau$, its stability at a branch point.  These points are
  indicated on the trivial solution branch with a dot for
  $\tau$ equal to
  $1/2, 1, 3/2$ and $2$ (left to right in the figure), where the latter is the model of Smith and
  co-authors.  The larger $\tau$, the more the bifurcation point
  shifts to the right. Thus the stability of the trivial solution
  is lost for larger values of the transport coefficient T. This
  result corresponds with the result in section
  \ref{subsec:StabilityOfTheTrivialSolution}. }

\change{The branches that emerge at these points contain steady state
  solutions with peaks as described earlier. Note that figure~\ref{fig:BifDiagram_tau} 
  only shows the bottom half of the
  emerging branches.  The branch at the right with $\tau=2$ contains
  the solutions with peaks for the model with exponent 2 and is the
  basic coupled model of Smith \textit{et al.}. This branch is the
  same branch as in figure \ref{fig:BifurcationDiagramScenario1}.
   For other choices of the exponents the figure is qualitatively
  the same. Only the stability changes slightly.  When $\tau$ is smaller
  the range of transport coefficients T where the solution branch is
  stable becomes smaller.}

\change{We have also found that the spacing between the peaks hardly
  changes if $\tau$ is changed in a continuous way}
 
\change{ So we conclude that there is little difference between the quadratic
  or linear dependence of the active transport on the auxin
  concentration.  Varying the exponent leads to qualitatively the same
  bifurcation.}

\begin{figure}
	\centering
	\includegraphics[height=0.9\textwidth,angle=-90]{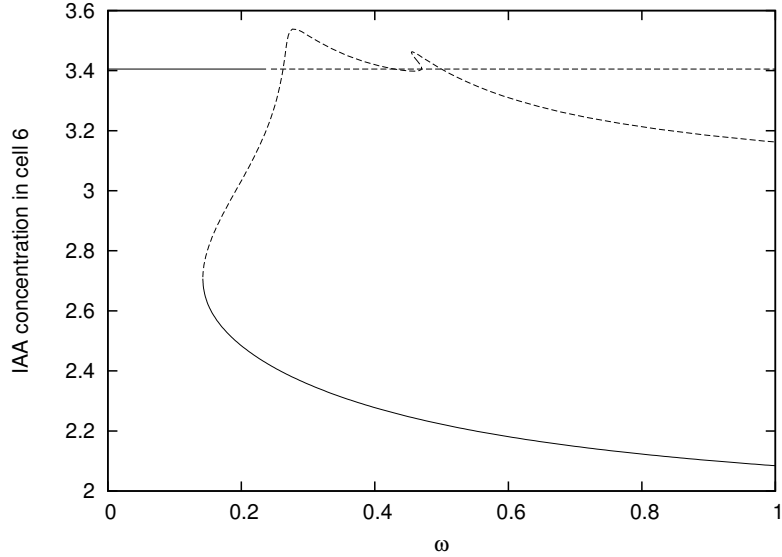}
  
        \caption{\change{The bifurcation diagram of example $6$ for a row of
          $20$ cells, parameter set M1 and active transport equation
          \eqref{eq:ActiveTransport_all} with $\tau = 2$. The diagram
          shows the IAA
          concentration in cell number $6$ versus the continuation parameter $\omega$. 
          The flat branch is the
          trivial solution, the other branch is found by starting with
          the steady state solution found in section
          \ref{subsec:TimeIntegration} for $\omega = 1$.}}
  \label{fig:BifDiagram_omega}
\end{figure}

\change{ 
  \paragraph{Example 6.}  This example explores the influence of the
  exponential dependence of the localization of PIN1 on the
  concentration of IAA.  We change, in a continuous way, the parameter
  $\omega$ in the active transport equation
  \eqref{eq:ActiveTransport_all}. Figure \ref{fig:BifDiagram_omega}
  shows the bifurcation diagram that depicts the concentration of
  auxin in cell number 6 versus the parameter $\omega$ for parameter
  set M1 and a row of $20$ cells. The transport
  coefficient T is now fixed at $1.5$.}
\change{ Again the flat branch represents the trivial solution that is
  independent of the parameter $\omega$. When $\omega$ is equal to
  zero, the active transport equation is equal to equation
  \eqref{eq:ActiveTransport3} and here the trivial solution is stable
  (see also in section  \ref{subsec:StabilityOfTheTrivialSolution}). When $\omega$ increases, 
  and the model transitions into the basic coupled
  model of Smith \textit{et al.}, the trivial solution becomes
  unstable (at $\omega = 0.2371$).  The other branch on this figure
  contains steady state solutions with peaks. For $\omega$ equal to
  one, the solution corresponds with the steady state solution in
  figure \ref{fig:timeEvolution}. When $\omega$ becomes smaller, the
  solution with peaks loses its stability in a turning point ($\omega
  = 0.1424$). The solution with peaks of figure
  \ref{fig:timeEvolution} in the model with $\omega = 1$ has no
  corresponding solution in the model with $\omega = 0$. This
  bifurcation diagram shows only a connection between the model with
  $\omega = 0$ and $\omega = 1$ by the trivial solution.}

\change{  Similar figures are found for other choices of the transport
  coefficient $T$.
}

%
%
%
%
\section{Conclusion and discussion}
\label{sec:ConclusionAndDiscussion}
In this paper we explored the model proposed by \citet{Smith2006}
for the transport of \change{auxin} in a one-dimensional row of
cells. The model describes the evolution of the
concentration of PIN1 and IAA in each of the cells by a coupled set of
non-linear ordinary differential equations.  The change in
concentration of IAA in each cell depends only on PIN1 and IAA
concentrations in that cell, its nearest and next-nearest neighbors
and this leads to a sparsely coupled system.

We analyzed the steady state solutions of the system as a
function of three of the 11 parameters in the problem.  We have varied
the IAA transport coefficient T, the diffusion coefficient D and the
production rate $\rho_{\IAA}$. \change{Furthermore, we introduced
  a slight generalization of the active transport equation.}  The main tool used
in the paper is numerical continuation to generate these solutions
starting from a trivial solution of the system.

The trivial solution is identified as an analytical solution where the
concentration in all the cells is the same. This solution is stable
for some region of the parameter space.  However, changing the
parameters, for example, increasing the IAA transport coefficient T,
destroys its stability.  In contrast to the uncoupled system studied
by \citet{Jonsson2006}, the eigenvalues of the Jacobian of
this coupled system can not be easily analyzed analytically. The
Jacobian now has a blocked sparse structure and its eigenvalues are
studied numerically.

In the exploration of the solutions, we identified two
generic bifurcation scenarios through which the trivial solution loses
its stability.  These scenarios reappear for various choices of the
parameters. In the first scenario, a stable solution can lose its
stability through a branch point, where it becomes a pattern with
regular spaced peaks of high auxin concentration.  This solution was
already found by Smith and collaborators by direct numerical
simulation.  The spacing and the height of the peaks in the pattern
depends on the other parameters of the system.

However, we have found in the coupled system that the trivial state
can also lose it stability through a Hopf bifurcation, where the
Jacobian has two complex conjugate eigenvalues that become purely imaginary.
For a limited parameter range this leads to stable periodic solutions,
where the concentration in the cell changes periodically over time.
However, for another range of parameters these periodic orbits are
unstable and the trivial solution then loses its stability through an
unstable orbit. The steady state solution then falls, beyond a
parameter threshold, back to a pattern of regularly spaced
peaks (see figure \ref{fig:bifurcationdiagramsolutionshopf}).  These
Hopf bifurcation and the periodic solutions are not present in the
model studied by J\"onsson where the PIN1 concentration is kept
constant and all eigenvalues of the Jacobian are real.

These are the only two bifurcation scenarios that we have found at the
stability boundary of the trivial solution for various choices of the
11 parameters in the model and for an increasing number of cells in the
row.

\change{We also explored modifications to the model of the active
  transport. The quadratic dependence on the auxin concentration is
  replaced by a linear dependence. The resulting bifurcation diagram
  is qualitatively the same. However, replacing the exponential
  dependence of the pin localization on the auxin transport leads to a
  different dynamical system where the trivial solution is always
  stable.}

Although the paper studies the steady state solutions of a rather
academic model with a row of equal sized square cells, the authors
believe it is a valuable contribution to our understanding of pattern formation by
auxin accumulation since it
builds the foundation for a rigorous bifurcation analysis of the
steady state patterns in a two dimensional array of cells.
 
\change{As demonstrated by \citet{Smith2006} generalized linear models 
are a powerful first step to understand the behavior of the biological system. Such studies
can be complemented by simulations of patterns of higher complexity that more closely resemble 
specific biological systems. }

\change{Applying the numerical continuation framework to 
a more general network of connected (ir)regular shaped cells, 
requires a Newton-Krylov solver so that the Newton
iteration deals with the non-linearity and the Krylov iteration 
solves the Jacobian system exploiting the sparsity. 
Furthermore the Krylov subspace iteration only requires the application 
of the Jacobian to a vector, which avoids  the explicit construction 
of the Jacobian matrix \citep{Kelley1995}.}

For the two dimensional array of cells, we expect to find a similar trivial solution that will lose again its
stability as the parameters change and turn into regular patterns of
high concentration peaks and time periodic solutions.

\change{In this paper we, essentially re-implemented the same linear file of cells as \citet{Smith2006},
but instead of using wrap-around boundary conditions, we implemented homogeneous Neumann boundary
conditions. In principle it is also possible to investigate the effect of varying in- and out-flux from 
the system and performing similar studies for these inhomogeneous Neumann boundary conditions. However, 
then there is no analytically solvable trivial solution. }

\change{There are several situations where changing boundary
  conditions are biologically relevant. First, in the leaf mesophyll,
  where a row of cells in which local maxima may induce new vascular
  development could be located between margins. These are a possible source of
  auxin and emerging veins that drain the auxin \citep{Scarpella2006},
  or at a later stage be surrounded by veins draining locally produced
  auxin. Second, in epidermal cells at the shoot apex along a tip to
  base gradient, where auxin comes in from the base and is transported
  to apical positions.  And finally, in central linear cell files in
  the root, auxin flows from base to tip and in the external layers
  auxin flows from the tip to the base leading to auxin gradients
  \citep{Grieneisen2007}.  }

We have explored continuation in
\change{the Neumann boundary condition} and found, amongst others,
s-shaped bifurcation diagrams with double limit points.  This leads to
a hysteresis effect in the boundary conditions. \change{This is
  described in \cite{Draelants2012}}.
  
\change{Also for problems with periodic boundary conditions a similar
  analysis can be performed. There, however, the Jacobian will have a
  one dimensional null-space associated with the translation
  invariance of solutions.  This can be regularized by introducing a
  phase condition as described by \citep{champneyssandstede}.}

There are still many uncertainties in the current generation of
models. Especially the large number of parameters and the uncertainty
in their values is a reason for concern.  By focusing on the
qualitative properties of the transitions that appear in the models
rather than on the states for particular choices of the parameters,
we hope to understand more about the possible patterns that appear in
real systems. It is valuable to calculate similar bifurcation
diagrams for all the proposed models for auxin transport using the
numerical continuation methods. \change{We could also use homotopy and 
use a continuous transformation to go from one model into another model.
Numerical continuation can then follow the solutions from one model
into the solutions of another model.} This will allow the comparison of
models across a range of parameters and check if they exhibit
qualitatively the same transition if parameters change.

In real plants it is impossible to tune a parameter such as the transport
coefficient T in a continuous way as is done in these calculations. It
can only be changed in discrete steps in a plant by the introduction
of, for example, mutations that compromise or enhance the auxin production or
transport.  Comparing such experiments with the model outcome will
make it possible to refine the models to give a more realistic
description of the biological system.

The bifurcation analysis \change{on the linear system with zero-influx at the boundaries}
yields interesting new insights into the
potential behaviors of the biological system.  It is interesting to
note that low values of the IAA transport coefficient lead to flat
distributions of the auxin concentration, whereas high concentrations
are required to establish sharp accumulation peaks. 

\change{ 1-N-naphthylphthalamic acid (NPA) has long been used
  experimentally as an inhibitor of auxin efflux. The precise
  molecular mechanism by which NPA impacts on auxin transport is still
  unknown, although several NPA-binding proteins have been identified
  \citep{muday}. Previously, we have assumed that NPA inhibits the
  cycling of PIN proteins to the cell membrane \citep{merks}. In terms
  of our model, this would relate to a lowering of $\rho_{_{\PIN_0}}$
  and/or $\rho_{_{\PIN}}$ which results in the formation of lower auxin peaks.
  Indeed, experimental inhibition of auxin transport with NPA abolishes the
normal narrow accumulation peaks and results in a much flatter auxin
distribution pattern \citep{Scarpella2006}. Also consistently with the model, 
the vam3 mutation that perturbs
PIN1 polarization, prevents auxin peaks forming and inhibits formation 
of higher order veins \citep{Shirakawa2009}.
More recently it was shown that a mutation in the auxin import carrier lax2, 
which inhibits auxin transport increases the numbers of vascular 
strand breaks \citep{Peret2012}, 
which may also be explained by the model prediction of lower auxin maxima, 
which could be unable to fully induce vascular development. Thus, the model 
behavior appears to be consistent with physiological treatments and mutations 
that affect specific model parameters.}

\change{The bifurcation analysis also yields interesting new insights into 
additional potential behaviors of the biological system.} 
One especially interesting behavior is the oscillation obtained with
a specific set of parameter values (Fig \ref{fig:TimeEvolutionPeriodic}).  This behavior
has to our knowledge never been observed in the context of leaf
development. However, in the root basal meristem, oscillating auxin
concentrations have been observed and related to the regular
induction of laterals along the growing axis of the root
\citep{DeSmet2007}.

One other characteristic unveiled by the bifurcation analysis is that
across the region where a pattern of auxin accumulation peaks are
generated, these peaks occur at very stable distances. This would
imply that in the case of vascular patterning, the initial distance
between veins is relatively stable, implying that observed differences
in vascular density in mature leaves \citep{Dhondt2011} largely result
from differences in subsequent development.  This is an example of an
new hypothesis generated by modeling a biological process that can be
experimentally validated and underlines the importance of systematic
exploration of biologically relevant parameter variations.

It is important to repeat that we have kept the plant geometry fixed
in the current model. It is an open question how the calculations can
be extended to include cells to undergo growth and division.

\section*{Acknowledgments}
\label{sec:Acknowledgments}
We acknowledge fruitful discussions with Dirk De Vos and Przemyslaw
Klosiewicz.  DD acknowledges financial support from the Department of
Mathematics and Computer Science of the University of Antwerp.  

This work is part of the
    Geconcerteerde Onderzoeksactie (G.O.A.)  research grant "A System
    Biology Approach of Leaf Morphogenesis" granted by the research
    council of the University of Antwerp.

%
%
%
%

\begin{thebibliography}{16}
\providecommand{\natexlab}[1]{#1}
\providecommand{\url}[1]{{#1}}
\providecommand{\urlprefix}{URL }
\expandafter\ifx\csname urlstyle\endcsname\relax
  \providecommand{\doi}[1]{DOI~\discretionary{}{}{}#1}\else
  \providecommand{\doi}{DOI~\discretionary{}{}{}\begingroup
  \urlstyle{rm}\Url}\fi
\providecommand{\eprint}[2][]{\url{#2}}

\bibitem[{Allgower and Georg(1994)}]{Allgower}
Allgower E, Georg K (1994) Numerical Path Following. Springer-Verlag Berlin

\bibitem[{Benkov\'a et~al(2003)Benkov\'a, Michniewicz, Sauer, Teichmann,
  Seifertov\'a, J\"urgens, and Friml}]{Benkova2003}
Benkov\'a E, Michniewicz M, Sauer M, Teichmann T, Seifertov\'a D, J\"urgens G,
  Friml J (2003) Local, efflux-dependent auxin gradients as a common module for
  plant organ formation. Cell 115:591--602

\bibitem[{Bilsborough et~al(2011)Bilsborough, Runions, Barkoulas, Jenkins,
  Hasson, Galinha, Laufs, Hay, Prusinkiewicz, and Tsiantis}]{Bilsborough2011}
Bilsborough G, Runions A, Barkoulas M, Jenkins H, Hasson A, Galinha C, Laufs P,
  Hay A, Prusinkiewicz P, Tsiantis M (2011) Model for the regulation of
  arabidopsis thaliana leaf margin development. Proceedings of the National
  Academy of Sciences 108:3424--3429

\bibitem[{Champneys and Sandstede (2007)}]{champneyssandstede}
  Champneys A R and Sandstede B (2007) Numerical computation of coherent
  structures.  In: Numerical Continuation Methods for Dynamical
  Systems (Edited by B Krauskopf, HM Osinga and J Galan-Vioque).
  Springer 331--358

\bibitem[{Clewley et~al(2007)Clewley, Sherwood, LaMar, and Guckenheimer}]{PyDS}
Clewley R, Sherwood W, LaMar M, Guckenheimer J (2007) Pydstool, a software
  environment for dynamical systems modeling. url http://pydstool.sourceforge.net

\bibitem[{De~Smet et~al(2007)De~Smet, Tetsumura, De~Rybel, Frey, Laplaze,
  Casimiro, Swarup, Naudts, Vanneste, Audenaert, Inz\'e, Bennet, and
  Beeckman}]{DeSmet2007}
De~Smet I, Tetsumura T, De~Rybel B, Frey N, Laplaze L, Casimiro I, Swarup R,
  Naudts M, Vanneste S, Audenaert D, Inz\'e D, Bennet M, Beeckman T (2007)
  Auxin-dependent regulation of lateral root positioning in the basal meristem
  of arabidopsis. Development 134:681--690

\bibitem[{Dhondt et~al(2011)Dhondt, Van~Haerenborgh, Van~Cauwenbergh, Merks,
  Philips, Beemster, and Inz\'e}]{Dhondt2011}
Dhondt S, Van~Haerenborgh D, Van~Cauwenbergh C, Merks R, Philips W, Beemster G,
  Inz\'e D (2011) Quantitative analysis of venation patterns of arabidopsis
  leaves by supervised image analysis. The plant Journal 69:553--563

\bibitem[{Doedel et~al(1997)Doedel, Champneys, Fairgrieve, Kuznetsov,
  Sandstede, and Wang}]{Doedel1997}
Doedel E, Champneys A, Fairgrieve T, Kuznetsov Y, Sandstede B, Wang X (1997)
  Continuation and bifurcation software for ordinary differential equations
  (with homcont). Available by anonymous ftp from ftp cs concordia ca,
  directory pub/doedel/auto
  
\bibitem[{Draelants et~al(2012)Draelants, Vanroose, Broeckhove, Beemster}]{Draelants2012}
Draelants D, Vanroose W, Broeckhove J, Beemster G T S,
Influence of an exogenous model parameter on the steady states in an auxin transport model. Submitted 2012.  

\bibitem[{Grieneisen et~al(2007)Grieneisen, Xu, Mar{\'e}e, Hogeweg, Scheres}]{Grieneisen2007}
Grieneisen V A, Xu J, Mar{\'e}e A F M, Hogeweg P, Scheres B (2007) Auxin transport is sufficient to
generate a maximum and gradient guiding root growth. Nature 449:1008--1013


\bibitem[{Hairer et~al(2009)Hairer, N{\o}rsett, and Wanner}]{Hairer2009}
Hairer E, N{\o}rsett S, Wanner G (2009) Solving ordinary differential equations
  I : nonstiff problems. Springer

\bibitem[{Hoyle(2006)}]{Hoyle} 
Hoyle R B , (2006) Pattern formation: an introduction to methods. Cambridge University Press.

\bibitem[{J{\"o}nsson et~al(2006)J{\"o}nsson, Heisler, Shapiro, Meyerowitz, and
  Mjolsness}]{Jonsson2006}
J{\"o}nsson H, Heisler M, Shapiro B, Meyerowitz E, Mjolsness E (2006) An
  auxin-driven polarized transport model for phyllotaxis. PNAS
  103(5):1633--1638

\bibitem[{Kelley(1995)}]{Kelley1995} Kelley, C.T. (1995) Iterative
  methods for linear and nonlinear equations, Society for Industrial
  Mathematics


\bibitem[{Krauskopf et~al(2007)Krauskopf, Osinga, and
  Gal{\'a}n-Vioque}]{Krauskopf}
Krauskopf B, Osinga H, Gal{\'a}n-Vioque J (2007) Numerical continuation methods
  for dynamical systems: path following and boundary value problems. Springer
  Verlag

\bibitem[{Merks~et~al(2007)Merks, Van de Peer, Inz\'e and
    Beemster}]{merks} Merks R M H, Van de Peer Y, Inz\'e D,
  Beemster G T S (2007) Canalization without flux sensors: a
  traveling-wave hypothesis. Trends in Plant Science 12, 384-390


\bibitem[{Muday and Delong(2001)Muday and Delong}]{muday} Muday
  G K and DeLong A (2001) Polar auxin transport: controlling where
  and how much. Trends Plant Sci. 6, 535–542

\bibitem[{Palme and G\"alweiler(1999)}]{Palme1999}
Palme K, G\"alweiler L (1999) Pin-pointing the molecular basis of auxin
  transport. Current Opinion in Plant Biology 2(5):375--381

\bibitem[{P{\'e}ret et~al(2012)P{\'e}ret, Swarup, Ferguson, Seth, Yang, Dhondt, James, Casimiro, Perry, Syed, Yang, Reemer, Venison, Howell, Perez-Amador, Yun, Alonso, Beemster, Laplaze, Murphy, Bennett, Nielsen, Swarup}]{Peret2012}
P{\'e}ret B, Swarup K, Ferguson A, Seth M, Yang Y, Dhondt S, James N, Casimiro I, Perry P, Syed A, Yang H, Reemer J, Venison E, Howell C,  Perez-Amador M A, Yun J, Alonso J, Beemster G T S, Laplaze L, Murphy A, Bennett M J, Nielsen E, Swarup R. AUX/LAX genes encode a family of auxin influx transporters that perform distinct function during Arabidopsis development. Submitted to Plant Cell in 2012.

\bibitem[{Reinhardt et~al(2003)Reinhardt, Pesce, Stieger, Mandel,
  Baltensperger, Bennett, Traas, Friml, and Kuhlemeier}]{Reinhardt2003}
Reinhardt D, Pesce E, Stieger P, Mandel T, Baltensperger K, Bennett M, Traas J,
  Friml J, Kuhlemeier C (2003) Regulation of phyllotaxis by polar auxin
  transport. Nature 426(6964):255--260

\bibitem[{Salinger et~al(2005)Salinger, Burroughs, Pawlowski, Phipps, and
  Romero}]{Salinger2005}
Salinger A, Burroughs E, Pawlowski R, Phipps E, Romero L (2005) Bifurcation
  tracking algorithms and software for large scale applications. International
  Journal of Bifurcation Chaos in Applied Sciences and Engineering
  15(3):1015--1032

\bibitem[{Scarpella et~al(2006)Scarpella, Marcos, Friml, and
  Berleth}]{Scarpella2006}
Scarpella E, Marcos D, Friml J, Berleth T (2006) Control of leaf vascular
  patterning by polar auxin transport. Genes \& Development 20:1015--1027

\bibitem[{Seydel(1994)}]{Seydel}
Seydel R (1994) Practical bifurcation and stability analysis: from equilibrium
  to chaos, vol~5. Springer

\bibitem[{Shirakawa et~al(2009)Shirakawa, Ueda, Shimada, Nishiyama, Hara-Nishimura}]{Shirakawa2009}
Shirakawa M, Ueda H, Shimada T, Nishiyama C, Hara-Nishimura I (2009) Vacuolar SNAREs function in the 
formation of the leaf vascular network by regulating auxin distribution. Plant \& Cell Physiology 50(7):1319--1328. 

\bibitem[{Smith et~al(2006)Smith, Guyomarc'h, Mandel, Reinhardt, Kuhlemeier,
  and Prusinkiewicz}]{Smith2006}
Smith R, Guyomarc'h S, Mandel T, Reinhardt D, Kuhlemeier C, Prusinkiewicz P
  (2006) A plausible model of phyllotaxis. PNAS 103(5):1301--1306



\end{thebibliography}

\end{document}